\documentclass[11pt]{amsart}
\usepackage{diagrams}

\usepackage{amsmath}
\usepackage{amsthm}
\usepackage{amssymb}
\usepackage{amsfonts}
\usepackage{amsxtra}
\usepackage{amscd}
\usepackage{epsfig}
\usepackage{verbatim}
\usepackage{latexsym,amstext,epsfig}

\newsymbol\pp 1275
\newcommand{\Hom}{\operatorname{Hom}}
\newcommand{\End}{\operatorname{End}}
\newcommand{\Ext}{\operatorname{Ext}}
\newcommand{\ext}{\operatorname{ext}}
\newcommand{\Rep}{\operatorname{Rep}}
\newcommand{\SI}{\operatorname{SI}}
\newcommand{\SL}{\operatorname{SL}}
\newcommand{\GL}{\operatorname{GL}}
\newcommand{\ZZ}{\mathbb Z}
\newcommand{\CC}{\mathbb C}
\newcommand{\RR}{\mathbb R}
\newcommand{\NN}{\mathbb N}

\newcommand{\Id}{\operatorname{Id}}
\newcommand{\mult}{\operatorname{mult}}
\newcommand{\supp}{\operatorname{supp}}
\newcommand{\Mat}{\operatorname{Mat}}

\newtheorem{theorem}{Theorem}[section]
\newtheorem{proposition}[theorem]{Proposition}
\newtheorem{corollary}[theorem]{Corollary}
\newtheorem{lemma}[theorem]{Lemma}

\theoremstyle{definition}
\newtheorem{definition}[theorem]{Definition}
\newtheorem{remark}[theorem]{Remark}

\newtheorem{example}[theorem]{Example}

\title[]{Quivers, Long exact sequences and Horn type inequalities}

\author{Calin Chindris}
\address{University of Minnesota, School of Mathematics, Minneapolis, MN, USA}
\email{chindris@math.umn.edu}

\date{\today}

\begin{document}
\bibliographystyle{plain}
 %\keywords{...}
\begin{abstract}
We give necessary and sufficient inequalities for the existence of
long exact sequences of $m$ finite abelian $p$-groups with fixed
isomorphy types. This problem is related to some generalized
Littlewood-Richardson coefficients that we define in this paper.
We also show how this problem is related to eigenvalues of
Hermitian matrices satisfying certain (in)equalities. When $m =
3,$ we recover the Horn type inequalities that solve the
saturation conjecture for Littlewood-Richardson coefficients and
Horn's conjecture.
\end{abstract}
\maketitle
\section{Introduction}\label{intro}
\subsection{Motivation}
Our main motivation in this paper goes back to the celebrated
conjecture of A.~Horn \cite{H} on the possible eigenvalues of a
sum of two Hermitian matrices. As explained in W.~Fulton's paper
\cite{F1}, there are problems in other areas of mathematics that
have the exact same solution as the eigenvalues of sums of two
Hermitian matrices problem. Two of them are the problem concerning
the existence of \emph{short} exact sequences of finite abelian
$p$-groups and that of the non-vanishing of the
Littlewood-Richardson coefficients. To state these problems, we
recall some definitions first. For every partition
$\lambda=(\lambda_1, \dots, \lambda_r)$ and a (fixed) prime number
$p$, one can construct a finite abelian $p$-group
$M_{\lambda}=\ZZ/p^{\lambda_1}\times \cdots \times
\ZZ/p^{\lambda_r}$. It is known that every finite abelian
$p$-group is isomorphic to $M_{\lambda}$ for a unique $\lambda$.
We will say that such a group is an abelian $p$-group of type
$\lambda$.

Let $V$ be a complex vector space of dimension $n$. If $\lambda=(
\lambda_1, \dots, \lambda_n)$ is a weakly decreasing sequence of
$n$ integers we denote by $S^{\lambda}(V)$ the irreducible
rational representation of $\GL(V)$ with highest weight $\lambda.$
Given three weakly decreasing sequences
$\lambda(1),~\lambda(2),~\lambda(3)$ of $n$ integers, we define
the Littlewood-Richardson coefficient
$c_{\lambda(1),\lambda(3)}^{\lambda(2)}$ to be the multiplicity of
$S^{\lambda(2)}(V)$ in $S^{\lambda(1)}(V) \otimes
S^{\lambda(3)}(V)$, i.e.
$$
c_{\lambda(1),\lambda(3)}^{\lambda(2)}=\dim_{\mathbb C}
\Hom_{\GL(V)}(S^{\lambda(2)}(V), S^{\lambda(1)}(V) \otimes
S^{\lambda(3)}(V)).
$$

An $n \times n$ complex matrix $H$ is said to be Hermitian if
$H=\overline{H}^t.$ It is a basic fact that all the eigenvalues of
a Hermitian  matrix are real numbers. We always write the
eigenvalues of a Hermitian matrix in weakly decreasing order.

Now we can state the three problems mentioned above.

\textbf{P1. Short exact sequences.} For which partitions
$\lambda(1),~\lambda(2)$ and $\lambda(3)$ with at most $n$ parts,
does there exist a short exact sequence
$$0\to M_1\to M_2\to M_3\to
0 ,$$  where $M_i$ is a finite abelian $p$-group of type
$\lambda(i)$ for every $1\leq i\leq 3$.

\textbf{P2. Littlewood-Richardson coefficients.} For which weakly
decreasing sequences $\lambda(1),~\lambda(2)$ and $\lambda(3)$ of
$n$ integers, do we have that
$$c_{\lambda(1),\lambda(3)}^{\lambda(2)}\neq 0.$$

\textbf{P3. Eigenvalues of a sum.} For which weakly decreasing
sequences $\lambda(1),~\lambda(2)$ and $\lambda(3)$ of $n$  real
numbers, do there exist $n\times n$ complex Hermitian matrices
$H(1),H(2)$ and $H(3)$ with eigenvalues $\lambda(1)$, $\lambda(2)$
and $\lambda(3)$ respectively and
$$
H(2)=H(1)+H(3).
$$

The equivalence of Problems \textbf{P1} and \textbf{P2} is due to
Klein \cite{Kle}. In \cite{H}, Horn conjectured that the set of
solutions to Problem \textbf{P3} consists of triples of $n$-tuples
of real numbers arranged in decreasing order satisfying certain
linear homogeneous inequalities. In fact, the following result has
been proved (we refer to the Notation paragraph from the end of
this section for basic definitions and notations).

\begin{theorem}[\textbf{Horn's conjecture}]\label{Hornconj.intro}
Let $\lambda(i)=(\lambda_1(i), \dots, \lambda_n(i)),$ $i \in \{1,$
$2,3 \}$ be three weakly decreasing sequences of $n$ real numbers.
Then the following are equivalent:
\begin{enumerate}
\renewcommand{\theenumi}{\arabic{enumi}}
\item there exist $n\times n$ complex Hermitian matrices
$H(1),H(2)$ and $H(3)$ with eigenvalues $\lambda(1)$, $\lambda(2)$
and $\lambda(3)$ respectively and
$$
H(2)=H(1)+H(3);
$$
\item the numbers $\lambda_j(i)$ satisfy
$$|\lambda(2)|=|\lambda(1)| + |\lambda(3)|,$$
together with
$$
\sum_{j \in I_2}\lambda_j(2) \leq \sum_{j\in I_1}\lambda_j(1) +
\sum_{j\in I_3}\lambda_j(3)
$$
for every triple $(I_1, I_2, I_3)$ of subsets of $\{1, \dots, n
\}$ of the same cardinality $r$ with $r<n$ and $c_{\lambda(I_1),
\lambda(I_3)}^{\lambda(I_2)} \neq 0.$\newline Assume that
$\lambda(i)$ are weakly decreasing sequences of $n$ integers. Then
$(1)$ and $(2)$ are equivalent to:

\item the Littlewood-Richardson coefficient $c_{\lambda(1),
\lambda(3)}^{\lambda(2)}$ is not zero.\newline Assume that
$\lambda(i)$ are partitions with at most $n$ parts. Then $(1)-(3)$
are equivalent to:

\item there exists a short exact sequence
$$0\to M_1\to M_2\to M_3\to
0 ,$$  where $M_i$ is a finite abelian $p$-group of type
$\lambda(i)$ for every $1\leq i\leq 3.$

\end{enumerate}
\end{theorem}

The first step in solving Horn's conjecture was taken by A.
~Klyachko \cite{Kl} who found necessary and sufficient linear
homogeneous inequalities for the eigenvalue problem. This set of
solutions to Problem \textbf{P3} forms a rational convex
polyhedral cone $\mathcal K(n,3)$ in $\mathbb R^{3n},$ known as
the Klyachko's cone. In the same paper, Klyachko made the
connection between his solution to the eigenvalue problem and the
Littlewood-Richardson coefficients. The next step was taken by A.
~Knutson and T. ~Tao \cite{KT} who proved what is now known as the
\emph{Saturation Conjecture} for the Littlewood-Richardson
coefficients. Their proof is based on some combinatorial gadgets
called honeycombs. H.~Derksen and J.~Weyman \cite{DW1} proved the
Saturation Conjecture in the more general context of quiver
theory. In a subsequent paper \cite{KTW}, A.~Knutson, T.~Tao and
C.~Woodward have described all the facets of the Klyachko's cone.
This way, they have obtained a minimal list of Horn type
inequalities defining the Klyachko's cone:

\begin{theorem}\cite{KTW}\label{minimalist.intro}
The Klyachko's cone $\mathcal K(n,3)$ consists of triples
\linebreak $(\lambda(1), \lambda(2), \lambda(3))$ of weakly
decreasing sequences of $n$ real numbers for which
$$
|\lambda(2)|=|\lambda(1)|+|\lambda(3)|
$$
and
$$
\sum_{j \in I_2}\lambda_j(2) \leq \sum_{j\in I_1}\lambda_j(1) +
\sum_{j\in I_3}\lambda_j(3)$$ for every triple $(I_1, I_2, I_3)$
of subsets of $\{1, \dots, n \}$ of the same cardinality $r$ with
$r<n$ and $c_{\lambda(I_1), \lambda(I_3)}^{\lambda(I_2)}=1;$
furthermore, this is now a minimal list.
\end{theorem}

As shown in \cite{CC4}, \cite{CC}, \cite{CB}, and \cite{DW2} most
of the above results proved by Klyachko, Knutson, Tao and Woodward
can be naturally obtained using quiver theory.

\subsection{The generalized problems}
When focusing on the existence of \emph{short} exact sequences, it
seems natural to extend Problem $\textbf{P1}$ to the case of
\emph{long} exact sequences with zeros at the ends of finite
abelian $p$-groups. Since a long exact sequence breaks into short
exact sequences, we will replace the Littlewood-Richardson
coefficient in Problem \textbf{P2} with a sum of products of
Littlewood-Richardson coefficients.

Let $m\geq 3$ and $n\geq 1$ be two integers.

\begin{definition} Given $m$ weakly decreasing sequences
$\lambda(1), \dots, \lambda(m)$ of $n$ integers, the
\emph{generalized Littlewood-Richardson coefficient}
$f(\lambda(1),\dots, \lambda(m))$ is defined as follows:
$$f(\lambda(1), \dots, \lambda(m))
=\sum c_{\lambda(1),\mu(1)}^{\lambda(2)}\cdot
c_{\mu(1),\mu(2)}^{\lambda(3)}\cdots
c_{\mu(m-4),\mu(m-3)}^{\lambda(m-2)}\cdot
c_{\mu(m-3),\lambda(m)}^{\lambda(m-1)},
$$
where the sum is taken over all partitions $\mu(1), \dots,
\mu(m-3)$ with at most $n$ parts.
\end{definition}

The convention is that when $m=3,$
$f(\lambda(1),\lambda(2),\lambda(3))$ is the Littlewood-Richardson
coefficient $c_{\lambda(1),\lambda(3)}^{\lambda(2)}.$

As it turns out, the generalized Littlewood-Richardson
coefficients are also related with parabolic affine
Kazhdan-Lusztig polynomials and decomposition numbers for
$q$-Schur algebras. This will be explained in Section
\ref{reptheoretic}.

% If $\lambda(1), \dots,
%\lambda(m)$ are partitions, let $n$ be a positive integer such
%that $\lambda(1), \dots, \lambda(m)$ have at most $n$ parts.
%Viewing $\lambda(1), \dots, \lambda(m)$ as weakly decreasing
%sequences of $n$ integers, $f(\lambda(1), \dots, \lambda(m))$ does
%not depend on $n.$ In fact, the sum defining $f(\lambda(1), \dots,
%\lambda(m))$ can be taken over all partitions $\mu(1),$ $\dots,$
%$\mu(m-3)$ with at most $l = \max\{ l(\lambda(i)) \mid 1 \leq i
%\leq m \}$ parts.

Now we are ready to state our generalized problems.

\textbf{Q1. Long exact sequences.} For which partitions
$\lambda(1), \dots, \lambda(m)$ with at most $n$ parts, does there
exist a long exact sequence
$$0\to M_1\to M_2\to \dots \to M_m\to
0 ,$$  where $M_i$ is a finite abelian $p$-group of type
$\lambda(i)$ for every $1\leq i\leq m$.

\textbf{Q2. Generalized Littlewood-Richardson coefficients.} For
which weakly decreasing sequences $\lambda(1), \dots, \lambda(m)$
of $n$ integers, do we have that
$$f(\lambda(1), \dots, \lambda(m))\neq 0.$$

\textbf{Q3. Generalized eigenvalue problems.} For which weakly
decreasing sequences $\lambda(1), \dots, \lambda(m)$ of $n$ real
numbers, do there exist $n\times n$ complex Hermitian matrices
$H(1), \dots, H(m)$ with eigenvalues $\lambda(1), \dots,
\lambda(m)$ and
$$\sum_{i~even}H(i)=\sum_{i~odd}H(i);$$
if $m > 3$ we also have that
$$\sum_{1\leq j \leq i}(-1)^{i+j}H(j)~~ \text{has non-negative
eigenvalues},$$ for every $2 \leq i \leq m-2.$

Note that what makes Problem \textbf{Q3} different from Problem
\textbf{P3} are the conditions on the eigenvalues of the
alternating partial sums obtained when $m>3.$

\subsection{Statement of the results}
Our first result is the following saturation property of the
generalized Littlewood-Richardson coefficients:

\begin{theorem}[\textbf{Saturation property}]\label{satgen}
Let $\lambda(1), \dots, \lambda(m)$ be $m$ weakly decreasing
sequences of $n$ integers. Then for every integer $r\geq 1$ we
have
$$
f(\lambda(1), \dots, \lambda(m))\neq 0 \Longleftrightarrow
f(r\lambda(1), \dots, r\lambda(m))\neq 0.
$$
\end{theorem}

Next, we relate the generalized Littlewood-Richardson coefficients
with the generalized spectral problem above.

\begin{definition}
Let $\mathcal K(n,m)\subseteq \mathbb R^{nm}$ be the solution set
to Problem \textbf{Q3,} i.e, $\mathcal K(n,m)$ is the set of all
$m$-tuples $(\lambda(1),$ $ \dots,$ $\lambda(m))$ of weakly
decreasing sequences of $n$ reals for which there exist $n\times
n$ complex Hermitian matrices $H(i),~ i \in \{1, \dots, m\}$
satisfying the conditions of Problem \textbf{Q3.} We call
$\mathcal K(n,m)$ the \emph{generalized Klyachko's cone}.
\end{definition}

To describe the generalized Klyachko's cone, we need to introduce
some notation. Let $(I_1, \dots, I_m)$ be an $m$-tuple of subsets
of $\{1, \dots, n \}$ such that at least one of them has
cardinality at most $n-1.$ We define the following weakly
decreasing sequences of integers (using conjugate partitions):
\newline

$\underline\lambda(I_1)=\lambda'(I_1),\hspace{15pt}
\underline\lambda(I_m)=
\begin{cases}
\lambda'(I_m) & \text{if $m$ is odd} \\
\lambda'(I_m\setminus \{n\}) & \text{if $m$ is even},
\end{cases}
$

and for $2\leq i\leq m-1$ \newline

$\underline\lambda(I_i)=
\begin{cases}
\lambda'(I_i) & \text{if $i$ is even} \\
\lambda'(I_i)-((|I_i|-|I_{i+1}|-|I_{i-1}|)^{n-|I_{i}|}) & \text{if
$i \leq m-2$ is odd}\\
\lambda'(I_i)-((|I_{m-1}|-|I_{m-2}|-|I_{m}\setminus
\{n\}|)^{n-|I_{i}|}) & \text{if $i = m-1$ is odd.}
\end{cases}
$
\newline

Now, we can state our \textbf{generalization of Horn's
conjecture:}

\begin{theorem}\label{mainthm}
Let $\lambda(i)=(\lambda_1(i),\dots,\lambda_n(i)),~i \in \{1,
\dots, m \}$ be $m$ weakly decreasing sequences of $n$ real
numbers. Then the following are equivalent:
\begin{enumerate}
\renewcommand{\theenumi}{\arabic{enumi}}

\item $(\lambda(1),$ $\dots,$ $\lambda(m)) \in \mathcal K(n,m);$

\item the numbers $\lambda_j(i)$ satisfy
$$
\sum_{i~even}|\lambda(i)|=\sum_{i~odd}|\lambda(i)|
$$
together with
$$
(*) \hspace{15pt}\sum_{i~even} \left(\sum_{j \in
I_i}\lambda_j(i)\right) \leq \sum_{i~odd} \left(\sum_{j\in
I_i}\lambda_j(i)\right)
$$
for every $m$-tuple $(I_1, \dots, I_m)$ for which $|I_1| = |I_2|,$
$|I_{m-1}| = |I_m|,$ $\underline{\lambda}(I_i),~ 1 \leq i \leq m$
are partitions and
$$f(\underline{\lambda}(I_1), \dots,
\underline{\lambda}(I_m))\neq 0.$$
Assume that $\lambda(i)$ are
sequences of integers. Then $(1)-(2)$ are equivalent to:

\item $f(\lambda(1), \dots, \lambda(m)) \neq 0.$ \newline Assume
that $\lambda(i)$ are partitions. Then $(1)-(3)$ are equivalent
to:

\item there exists a long exact sequence of the form
$$0\to M_1\to M_2\to \dots \to M_m\to
0 ,$$  where $M_i$ is a finite abelian $p$-group of type
$\lambda(i)$ for every $1\leq i\leq m.$
\end{enumerate}
\end{theorem}

Note that the above Theorem gives a recursive method for finding
all non-zero generalized Littlewood-Richardson coefficients. It
turns out that one can shorten the list of Horn type inequalities
of Theorem \ref{mainthm}{(2)}:

\begin{proposition} \label{genKlycone}
The following statements are true.
\begin{enumerate}
\renewcommand{\theenumi}{\arabic{enumi}}
\item We have $$\dim \mathcal K (n, m)=mn-1.$$

\item The cone $\mathcal K (n, m)$ consists of all $m$-tuples
$(\lambda(1)$ $,\dots,$ $\lambda(m))$ of weakly decreasing
sequences of $n$ reals for which  $$
\sum_{i~even}|\lambda(i)|=\sum_{i~odd}|\lambda(i)|
$$ and $(*)$ holds for
every $m$-tuple $(I_1,\dots, I_m)$ for which $|I_1| = |I_2|,$
$|I_{m-1}| = |I_m|,$ $\underline{\lambda}(I_i),~1 \leq i \leq m$
are partitions and
$$f(\underline{\lambda}(I_1), \dots,
\underline{\lambda}(I_m)) = 1.$$
\end{enumerate}
\end{proposition}

We want to point out that our results do not depend on the work of
Klyachko, Knutson and Tao. In fact, our strategy is to show first
that the non-vanishing of the generalized Littlewood-Richardson
coefficients is equivalent to the existence of non-zero
semi-invariants for the generalized flag quiver setting. Once we
have switched to quiver invariant theory, our main tool is a nice
description of the facets of the cone of effective weights for
quivers without oriented cycles which was proved by Derksen and
Weyman \cite{DW2}.

The paper is organized as follows. In Section \ref{prelim}, we
recall a certain saturation property for effective weights for
quivers which is due to Derksen and Weyman \cite{DW1}. The
generalized flag quiver setting is defined in Section
\ref{flagquiver} where we also prove the saturation property for
the generalized Littlewood-Richardson coefficients. A more
detailed description of the so called cone of effective weights
for arbitrary quivers (without oriented cycles) is given in
Section \ref{coneabstract}. In Section \ref{facetsgenflag}, we
find the facets of the cone of effective weights associated to the
generalized flag quiver setting. The Horn type inequalities and
the $m$-tuples $(I_1, \dots, I_m)$ occurring in Theorem
\ref{mainthm}{(2)} are obtained in Section \ref{proofs}. In
Section \ref{proofsnou}, we give a moment map description of the
cone associated to the generalized flag quiver setting and  prove
Theorem \ref{mainthm} and Proposition \ref{genKlycone}. In Section
\ref{reptheoretic}, we discuss two representation theoretic
interpretations of the generalized Littlewood-Richardson
coefficients. First, we explain how the generalized
Littlewood-Richardson coefficients are related to some parabolic
affine Kazhdan-Lusztig polynomials and decomposition numbers for
$q$-Schur algebras. We also show how our coefficients can be
viewed as multiplicities of irreducible representations of a
product of general linear groups. In Section \ref{finalrmks}, we
make some comments on the minimality of our list of Horn type
inequalities.

\textbf{Notation.} A partition $\lambda$ of length $N$ is a
sequence of $N$ positive integers $\lambda=(\lambda_1, \dots,
\lambda_N)$ with $\lambda_1\geq \dots \geq \lambda_N \geq 1.$ We
say that $\lambda$ is a partition with at most $N$ (non-zero)
parts if $\lambda=(\lambda_1, \dots, \lambda_N)\in \ZZ^N$ with
$\lambda_1\geq \dots \geq \lambda_N \geq 0.$ A partition $\lambda$
will be also viewed as a weakly decreasing sequence of $n$
integers by adding zero parts, for any integer $n$ greater or
equal than the number of non-zero parts of $\lambda.$ If
$\lambda=(\lambda_1, \dots, \lambda_N)$ is a weakly decreasing
sequence then we define $r\lambda$ by $r\lambda=(r\lambda_1,
\dots, r\lambda_N).$ Let $\lambda = (\lambda_1, \dots, \lambda_N)$
and $\mu = (\mu_1, \dots, \mu_M)$ be two weakly decreasing
sequences of integers. Then we define the sum $\lambda + \mu$ by
first extending $\lambda$ or $\mu$ with zero parts (if necessary)
and then we add them componentwise. For a partition $\lambda,$ we
denote by $\lambda'$ the partition conjugate to $\lambda,$ i.e.,
the Young diagram of $\lambda'$ is the Young diagram of $\lambda$
reflected with respect to its main diagonal. We will often refer
to partitions as Young diagrams. If $I=\{z_1< \dots <z_r\}$ is an
$r$-tuple of integers then $\lambda(I)$ is defined by
$\lambda(I)=(z_r-r, \dots, z_1-1).$ For $r \geq 0$ and $a$ two
integers, we denote the $r$-tuple $(a, \dots, a)$ by $(a^r).$ If
$\lambda=(\lambda_1, \dots, \lambda_N)$ is a sequence of real
numbers, we define $|\lambda|=\sum_{i=1}^N \lambda_i.$

\section{Preliminaries} \label{prelim}
\subsection{Generalities}A quiver $Q=(Q_0,Q_1,t,h)$
consists of a finite set of vertices $Q_0$, a finite set of arrows
$Q_1$ and two functions $t,h:Q_1 \to Q_0$ that assign to each
arrow $a$ its tail $ta$ and its head $ha,$ respectively. We write
$ta{\buildrel a\over\longrightarrow}ha$ for each arrow $a \in
Q_1$.

For simplicity, we will be working over the field of complex
numbers $\CC.$ A representations $V$ of $Q$ over $\CC$ is a family
of finite dimensional $\CC$-vector spaces $\lbrace V(x) \mid x\in
Q_0\rbrace$ together with a family $\{ V(a):V(ta)\rightarrow V(ha)
\mid a \in Q_1 \}$ of $\CC$-linear maps. If $V$ is a
representation of $Q$, we define its dimension vector $\underline
d_V$ by $\underline d_V(x)=\dim_{\CC} V(x)$ for every $x\in Q_0$.
Thus the dimension vectors of representations of $Q$ lie in
$\Gamma=\ZZ^{Q_0}$, the set of all integer-valued functions on
$Q_0$. For each vertex $x$, we denote by $\varepsilon_x$ the
simple dimension vector corresponding to $x$, i.e.
$\varepsilon_x(y)=\delta_{x,y}~,~\forall~ y \in Q_0,$ where
$\delta_{x,y}$ is the Kronecker symbol.

Given two representations $V$ and $W$ of $Q$, we define a morphism
$\phi:V \rightarrow W$ to be a collection of linear maps $\lbrace
\phi(x):V(x)\rightarrow W(x)\mid x \in Q_0 \rbrace$ such that for
every arrow $a\in Q_1$, we have $\phi(ha)V(a)=W(a)\phi(ta)$. We
denote by $\Hom_Q(V,W)$ the $\CC$-vector space of all morphisms
from $V$ to $W$. In this way, we obtain the abelian category
$\Rep(Q)$ of all quiver representations of $Q.$ Let $W$ and $V$ be
two representations of $Q.$ We say that $V$ is a subrepresentation
of $W$ if $V(x)$ is a subspace of $W(x)$ for all vertices $x \in
Q_0$ and $V(a)$ is the restriction of $W(a)$ to $V(ta)$ for all
arrows $a \in Q_1.$

If $\alpha,\beta$ are two elements of $\Gamma$, we define the
Euler form by
\begin{equation}
\langle\alpha,\beta \rangle = \sum_{x \in Q_0}
\alpha(x)\beta(x)-\sum_{a \in Q_1} \alpha(ta)\beta(ha).
\end{equation}

\emph{From now on, we will assume that our quivers are without
oriented cycles.}

\subsection{Semi-invariants for quivers}
Let $\beta$ be a dimension vector of $Q$. The representation space
of $\beta-$dimensional representations of $Q$ is defined by
$$\Rep(Q,\beta)=\bigoplus_{a\in Q_1}\Hom(\CC^{\beta(ta)}, \CC^{\beta(ha)}).$$
If $\GL(\beta)=\prod_{x\in Q_0}\GL(\beta(x))$ then $\GL(\beta)$
acts algebraically on $\Rep(Q,\beta)$ by simultaneous conjugation,
i.e., for $g=(g(x))_{x\in Q_0}\in \GL(\beta)$ and $V=\{V(a)\}_{a
\in Q_1} \in \Rep(Q,\beta),$ we define $g \cdot V$ by
$$(g\cdot V)(a)=g(ha)V(a)g(ta)^{-1}\ \text{for each}\ a \in Q_1.$$ In
this way, $\Rep(Q,\beta)$ is a rational representation of the
linearly reductive group $\GL(\beta)$ and the $\GL(\beta)-$orbits
in $\Rep(Q,\beta)$ are in one-to-one correspondence with the
isomorphism classes of $\beta-$dimensional representations of $Q.$
As $Q$ is a quiver without oriented cycles, one can show that
there is only one closed $\GL(\beta)-$orbit in $\Rep(Q,\beta)$ and
hence the invariant ring $\text{I}(Q,\beta)= \CC
[\Rep(Q,\beta)]^{\GL(\beta)}$ is exactly the base field $\CC.$
Although there are only constant $\GL(\beta)-$invariant polynomial
functions on $\Rep(Q,\beta)$, the action of $\SL(\beta)$ on
$\Rep(Q,\beta)$ provides us with a highly non-trivial ring of
semi-invariants.

Let $\SI(Q,\beta)= \CC [\Rep(Q,\beta)]^{\SL(\beta)}$ be the ring
of semi-invariants. As $\SL(\beta)$ is the commutator subgroup of
$\GL(\beta)$ and $\GL(\beta)$ is linearly reductive, we have that
$$\SI(Q,\beta)=\bigoplus_{\sigma
\in X^\star(\GL(\beta))}\SI(Q,\beta)_{\sigma},
$$
where $X^\star(\GL(\beta))$ is the group of rational characters of
$\GL(\beta)$ and $$\SI(Q,\beta)_{\sigma}=\lbrace f \in \CC
[\Rep(Q,\beta)] \mid gf= \sigma(g)f, \text{for all}~g \in
\GL(\beta)\rbrace$$ is the space of semi-invariants of weight
$\sigma.$ Note that a character or weight of $\GL(\beta)$ is of
the form
$$\{g(x) \mid x \in Q_0\} \in \GL(\beta) \mapsto \prod_{x \in
Q_0}(\det g(x))^{\sigma(x)}$$ with $\sigma(x) \in \ZZ$ for all $x
\in Q_0.$ Therefore, we can identify $X^\star(\GL(\beta))$ with
$\ZZ ^{Q_0}.$ If $\alpha \in \Gamma$, we define $\sigma=\langle
\alpha,\cdot \rangle$ by
$$\sigma(x)=\langle \alpha,\varepsilon_x \rangle~,~\forall x\in
Q_0.$$ Similarly, one can define $\sigma = \langle \cdot,\alpha
\rangle.$

Given a quiver $Q$ and a dimension vector $\beta$, we define the
set $\Sigma(Q, \beta)$ of \emph{(integral) effective weights} by
$$\Sigma(Q,\beta)=\lbrace \sigma \in \ZZ^{Q_0} \mid
\SI(Q,\beta)_{\sigma}\neq 0\rbrace.$$

In \cite{S2}, Schofield constructed some very useful
semi-invariants for quivers. A fundamental result due to Derksen
and Weyman \cite{DW1} (see also \cite{SVB}) states that these
semi-invariants span all spaces of semi-invariants. An important
consequence of this spanning theorem is the following saturation
property.

\begin{proposition} \cite[Theorem 3]{DW1}\label{saturation}
If $Q$ is a quiver and $\beta$ is a dimension vector, then the set
$$\Sigma(Q,\beta)=\lbrace \sigma \in \ZZ^{Q_0} \mid
\SI(Q,\beta)_{\sigma}\neq 0\rbrace,$$ is saturated, i.e., if
$\sigma$ is a weight and $r \geq 1$  $$\SI(Q,\beta)_{\sigma}\neq 0
\Longleftrightarrow \SI(Q,\beta)_{r \sigma}\neq 0.$$
\end{proposition}

A detailed description of the set $\Sigma(Q,\beta)$ of effective
weights can be found in Section \ref{coneabstract}, Theorem
\ref{semi1} and Proposition \ref{descfacets}.

\section{The generalized flag quiver and the saturation
property}\label{flagquiver} In this section we first define the
generalized flag quiver and show that the generalized
Littlewood-Richardson coefficients are the dimensions of the
spaces of semi-invariants for the generalized flag quiver.

Let $m\geq 3$ and $n\geq 1$ be two positive integers. The
\emph{generalized flag quiver setting} is defined as follows.
\begin{enumerate}
\renewcommand{\theenumi}{\alph{enumi}}
\item The quiver $Q$ has $m-2$ \emph{central vertices} $2=(n,
2)=(n, 1), 3=(n,3), \dots, m-2=(n, m-2), m-1=(n, m-1)=(n, m)$ at
which we attach $m$ equioriented $\mathbb A_n$ quivers (or flags)
$\mathcal F(1), \dots, \mathcal F(m)$ such that $\mathcal F(i)$
goes in the corresponding cental vertex $(n, i)$ if $i$ is even
and it goes out from the corresponding cental vertex $(n, i)$ if
$i$ is odd. Furthermore, there are $m-3$ main arrows $a_1, \dots,
a_{m-3}$ connecting the central vertices such that $i+1{\buildrel
a_i\over\longrightarrow}i+2$ if $i$ is odd and $i+2{\buildrel
a_i\over\longrightarrow}i+1$ if $i$ is even. For example, if the
number of flags $m$ is even then our quiver $Q$ looks like
%\vfill\pagebreak
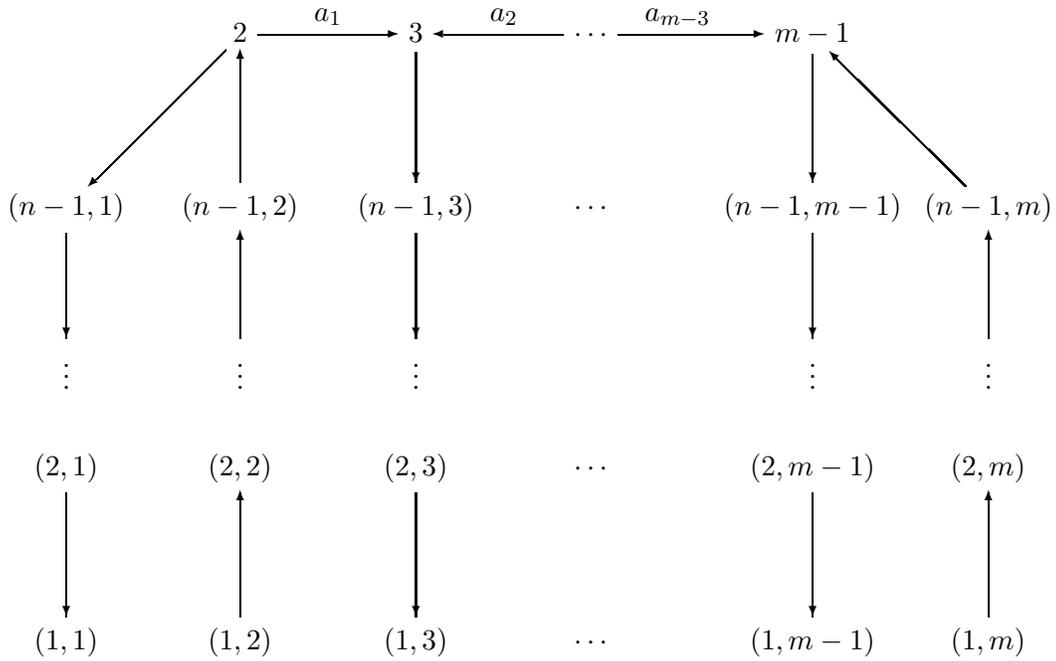
\begin{diagram}
& & 2& \rTo{a_1} &  3&\lTo{a_2}  & \cdots  &\rTo{a_{m-3}} & m-1 \\
& \ldTo & \uTo & & \dTo & & & &\dTo &\luTo\\
(n-1, 1)& & (n-1, 2)& & (n-1, 3)& &\cdots & & (n-1, m-1)&
& (n-1, m)\\
\dTo& &\uTo& &\dTo& & &&\dTo& &\uTo\\
\vdots& &\vdots& &\vdots& & & &\vdots& &\vdots\\
(2, 1)& & (2, 2)& & (2, 3)& &\cdots & & (2, m-1)& & (2, m)\\
\dTo& &\uTo& &\dTo& & & &\dTo& &\uTo\\
(1, 1)& & (1, 2)& & (1, 3)& &\cdots & & (1, m-1)& & (1, m)
\end{diagram}

\vfill\pagebreak

\item The dimension vector $\beta$ is given by $\beta(j, i)=j$ for
all $j \in \{1, \dots, n\}$ and $i \in \{1, \dots, m\},$ i.e.,
$\beta$ is equal to
$$
\begin{matrix}
& n & n & \cdots & n \\
n-1& n-1 & n-1 & \cdots & n-1& n-1\\
\vdots& \vdots& \vdots &  & \vdots& \vdots\\
2& 2 & 2 & \cdots & 2 &2\\
1& 1 & 1 & \cdots & 1 & 1
\end{matrix}$$
\end{enumerate}

In this section, the only quiver setting we will be working with
is the generalized flag quiver setting.

\begin{lemma}\label{compute} Let $\sigma \in \ZZ^{Q_0}$ be a weight. If $\dim\SI(Q,\beta)_{\sigma}\neq 0$
then:
\begin{enumerate}
\renewcommand{\theenumi}{\arabic{enumi}}
\item the weight $\sigma$ must satisfy the inequalities
$$(-1)^i\sigma(j,i) \geq 0,$$ for
all $1 \leq j \leq n,~2\leq i \leq m-1$ and
$$(-1)^i \sigma(j, i) \geq 0,$$ for all $1 \leq j \leq n-1,~i
\in \{1, m \};$

\item we have
$$\dim
\SI(Q,\beta)_{\sigma}=\sum_{\mu(1),\cdots,\mu(m-3)}c_{\gamma(1),\mu(1)}^{\gamma(2)}\cdot
c_{\mu(1),\mu(2)}^{\gamma(3)}\cdots
c_{\mu(m-3),\gamma(m-1)}^{\gamma(m)},$$ where
$$
\begin{aligned}
\gamma(1) &
=((n-1)^{-\sigma(n-1, 1)},\dots,1^{-\sigma(1,1)})', \\
\gamma(m) & =((n-1)^{(-1)^m \cdot
\sigma(n-1, m)},\dots,1^{(-1)^m \cdot \sigma(1, m)})', \\
\gamma(i) & =(n^{(-1)^i \cdot \sigma(n, i)},\dots,1^{(-1)^i \cdot
\sigma(1, i)})',
\end{aligned}
$$
for all $i \in \{2, \dots, m-1\}.$
\end{enumerate}
\end{lemma}

\begin{proof}
The first part of this Lemma follows as we compute
$\SI(Q,\beta)_{\sigma}.$ For simplicity, let us define $V_j(i)=
\CC^{\beta(j, i)}$. Using Cauchy's formula \cite[page 121]{F3}, we
can decompose the affine coordinate ring $\CC [\Rep(Q,\beta)]$ as
a sum of tensor products of irreducible representations of the
general linear groups $\GL(V_j(i))$. The idea is to identify those
terms that will give us non-zero semi-invariants of weight
$\sigma.$ An arbitrary term in this decomposition is made up of
tensor products of irreducible representations coming from the $m$
flags. If $\mathcal F(i)$ is a flag going in the cental vertex
$(n,i),$ then the $n-1$ arrows of this flag contribute with
$$
S^{\gamma^1(i)}V_1(i)\otimes \bigotimes_{j=2}^{n-1}
\left(S^{\gamma^{j-1}(i)}V^{*}_j(i)\otimes S^{\gamma^j(i)}V_j(i)
\right)\otimes S^{\gamma^{n-1}(i)}V^{*}_n(i),
$$
for partitions $\gamma^1(i), \dots, \gamma^{n-1}(i).$

When computing semi-invariants, we see that $\left
(S^{\gamma^1(i)}V_1(i)\right )^{\SL(V_1(i))}$ is non-zero if and
only if it is one dimensional. In this case, $\gamma^1(i)$ is a
$\beta(1, i)\times w$ rectangle and the space is spanned by a
semi-invariant of weight $w.$ So, $\left
(S^{\gamma^1(i)}V_1(i)\right )^{\SL(V_1(i))}$ contains non-zero
semi-invariants of weight $\sigma(1, i)$ if and only if $\sigma(1,
i) \geq 0$ and $\gamma^1(i) = (\sigma(1, i)^{\beta(1, i)}),$ i.e.,
$$\gamma^1(i) =
(1^{\sigma(1, i)})'.
$$
Next, we look at the space
$$ \left(S^{\gamma^1(i)}V^{*}_2(i)\otimes
S^{\gamma^2(i)}V_2(i) \right)^{\SL(V_2(i))}$$ which is canonically
isomorphic to
$\Hom_{\SL(V_2(i))}(S^{\gamma^1(i)}V_2(i),S^{\gamma^2(i)}V_2(i)).
$ Now, this space is non-zero if and only if it is one dimensional
in which case $\gamma^2(i)$ is $\gamma^1(i)$ plus some extra
columns of length $\beta(2, i)$ and the number of these extra
columns is the weight of a semi-invariant spanning this space.
Consequently, $\left(S^{\gamma^1(i)}V^{*}_2(i)\otimes
S^{\gamma^2(i)}V_2(i) \right)^{\SL(V_2(i))}$ contains non-zero
semi-invariants of weight $\sigma(2, i)$ if and only if $\sigma(2,
i)\geq 0$ and $\gamma^2(i)$ is $\gamma^1(i)$ plus $\sigma(2, i)$
columns of length $\beta(2, i),$ i.e.,
$$\gamma^2(i) = (2^{\sigma(2, i)},
1^{\sigma(1, i)})'.$$

Reasoning in this way, we see that the vertices of this flag
$\mathcal F(i)$, except the central one $(n, i)$, give non-zero
spaces of semi-invariants (in which case they must be one
dimensional) of weight $\sigma(1, i), \dots, \sigma(n-1, i)$ if
and only if $\sigma(j,i)\geq 0$ for all $1 \leq j \leq n-1,$
$\gamma^1(i)$ is a $\beta(1, i)\times \sigma(1, i)$ rectangle and
$\gamma^j(i)$ is $\gamma^{j-1}(i)$ plus $\sigma(j, i)$ columns of
length $\beta(j, i)$ for all $j \in \{2, \dots, n-1 \},$ i.e.,
$$
\gamma^{n-1}(i) = ((n-1)^{\sigma(n-1, i)},\dots,1^{\sigma(1,
i)})'.
$$
We have proved that a flag $\mathcal F(i)$ going in the central
vertex $(n,i)$ contributes to the space of semi-invariants $\SI(Q,
\beta)_{\sigma}$ with
$$
S^{\gamma^{n-1}(i)}V^{*}_n(i),
$$
where $\gamma^{n-1}(i)$ is completely determined by the weight
$\sigma$ along the flag $\mathcal F(i).$

Similarly, if $\mathcal F(l)$ is a flag going out from the central
vertex $(n,l),$ then $\sigma(j,l)\leq 0$ for all $1 \leq j \leq
n-1$ and $\mathcal F(l)$ contributes to $\SI(Q, \beta)_{\sigma}$
with
$$
S^{\gamma^{n-1}(l)}V_n(l),
$$
where
$$
\gamma^{n-1}(l) = ((n-1)^{-\sigma(n-1,l)},\dots,1^{-\sigma(1,
l)})'.
$$

Next, the main $m-3$ arrows of our quiver give us partitions
$\mu(1),\dots,\mu(m-3),$ with at most $n$ parts, and the central
vertices give us the following spaces of semi-invariants:
$$\left(S^{\gamma^{n-1}(1)}V(2)\otimes S^{\mu(1)}V(2)\otimes
S^{\gamma^{n-1}(2)}V^{*}(2) \right)^{\SL(V(2))}$$ coming from the
vertex $2,$ $$\left( S^{\gamma^{n-1}(3)}V(3) \otimes
S^{\mu(1)}V^{*}(3)\otimes S^{\mu(2)}V^{*}(3) \right)^{\SL(V(3))}$$
coming from the vertex $3$ and so on. Taking into account the
weights at the central vertices, it is clear that the dimension of
the space of semi-invariants $\SI(Q,\beta)_{\sigma}$ is the
desired sum of products of Littlewood-Richardson coefficients.
\end{proof}

Let $\lambda(1), \dots, \lambda(m)$ be weakly decreasing sequences
of $n$ integers. To show that $f(\lambda(1), \dots, \lambda(m))$
can be viewed as the dimension of a space of semi-invariants, we
are going to apply Lemma \ref{compute}. Let us define
$\sigma_{\lambda}$ by

\begin{equation} \label{thewt1}
\sigma_{\lambda}(j, i) = (-1)^i(\lambda_j(i)-\lambda_{j+1}(i)),
\forall~ 1 \leq j \leq n-1, \forall~1 \leq i \leq m,
\end{equation}

\begin{equation} \label{thewt2}
\sigma_{\lambda}(i) = (-1)^i\lambda_n(i), \forall~ i \neq 2,~ m-1,
\end{equation}

\begin{equation} \label{thewt3}
\sigma_{\lambda}(2) = \lambda_n(2)-\lambda_n(1),
\end{equation}

\begin{equation} \label{thewt4}
\sigma_{\lambda}(m-1)=(-1)^{m-1}(\lambda_n(m-1)-\lambda_n(m)).
\end{equation}

If $m=3$ then $\sigma_{\lambda}$ at the central vertex becomes
$$
\sigma_{\lambda}(2) = \lambda_n(2)-\lambda_n(1) - \lambda_n(3).
$$ With these notations we have:

\begin{lemma} \label{wtsigmalambda}
Let $\lambda(1), \dots, \lambda(m)$ be $m \geq 3$ weakly
decreasing sequences of $n$ integers. Then for every integer $r
\geq 1,$ we have
$$
f(r\lambda(1), \dots, r\lambda(m))=\dim \SI(Q,
\beta)_{r\sigma_{\lambda}}.
$$
\end{lemma}

\begin{proof}
We prove this Lemma when $r = 1,$ as the general case reduces to
this one. First, let us consider the following transformations
$$
\begin{aligned}
&\gamma(1)=\lambda(1)-(\lambda_n(1)^n),\\
&\gamma(2)=\lambda(2)-(\lambda_n(1)^n),\\
&\gamma(m-1)=\lambda(m-1)-(\lambda_n(m)^n),\\
&\gamma(m)=\lambda(m)-(\lambda_n(m)^n),\\
&\gamma(i)=\lambda(i), \forall ~i\notin \{1,2,m-1,m\}.
\end{aligned}
$$
If $m=3$ then $\gamma(2)$ becomes
$\gamma(2)=\lambda(2)-((\lambda_n(1)+\lambda_n(3))^n).$ With this
transformations, we have
$$
\begin{aligned}
\gamma(1) &
=((n-1)^{-\sigma(n-1, 1)},\dots,1^{-\sigma(1,1)})', \\
\gamma(m) & =((n-1)^{(-1)^m \cdot
\sigma(n-1, m)},\dots,1^{(-1)^m \cdot \sigma(1, m)})', \\
\gamma(i) & =(n^{(-1)^i \cdot \sigma(n, i)},\dots,1^{(-1)^i \cdot
\sigma(1, i)})',
\end{aligned}
$$
for all $i \in \{2, \dots, m-1\}.$ Applying Lemma \ref{compute},
we get that
$$
f(\gamma(1), \dots, \gamma(m))=\dim \SI(Q,
\beta)_{\sigma_{\lambda}}.
$$
On the other hand, we clearly have $ f(\lambda(1), \dots,
\lambda(m)) = f(\gamma(1), \dots, \gamma(m))$ and so the proof
follows.
\end{proof}

\begin{remark}
Let us note that if $f(\lambda(1), \dots, \lambda(m))$ is non-zero
then the first part of Lemma \ref{compute} tells us that
$\lambda(i), i \notin \{1,2,m-1,m \}$ are in fact partitions. Of
course, this is also clear from the definition of $f(\lambda(1),
\dots, \lambda(m)).$
\end{remark}

\begin{proof}[Proof of Theorem \ref{satgen}]
The proof follows from Proposition \ref{saturation} and Lemma
\ref{wtsigmalambda}.
\end{proof}

\section{The cone of effective weights for quivers} \label{coneabstract}

Let $Q$ be a quiver without oriented cycles and let $\beta$ be a
dimension vector. In this section we will further describe the
rational convex polyhedral cone whose lattice points form the set
of integral effective weights
$$\Sigma(Q,\beta)=\lbrace \sigma
\in \ZZ^{Q_0} \mid \SI(Q,\beta)_{\sigma}\neq 0\rbrace.$$

If $\sigma \in \mathbb R^{Q_0}$ is a real valued function on the
set of vertices and $\alpha \in \Gamma,$ we define
$\sigma(\alpha)$ by
$$\sigma(\alpha) = \sum_{x \in Q_0} \sigma(x)\alpha(x).$$

A necessary condition for a weight $\sigma \in \ZZ^{Q_0}$ to
belong to $\Sigma(Q,\beta)$ is $\sigma(\beta)=0.$ Indeed, the
action of the one dimensional torus $\{(t \Id_{\beta(i)})_{i \in
Q_0} \mid t \in K\setminus \{0\} \}$ on the representation space
$\Rep(Q,\beta)$ is trivial. If $f$ is a non-zero semi-invariant of
weight $\sigma$ and $g_t = (t \Id_{\beta(i)})_{i \in Q_0} \in
\GL(\beta)$ then
$$
g_t \cdot f = t^{\sigma(\beta)} \cdot f
$$
clearly implies that $\sigma(\beta) = 0.$
\begin{lemma}[Reciprocity Property]\cite[Corollary 1]{DW1}\label{reciprocity} We have
$$\dim\SI(Q,\beta)_{\langle \alpha, \cdot \rangle} = \dim
\SI(Q,\alpha)_{-\langle \cdot, \beta \rangle}.$$
\end{lemma}

In this case, we define $\alpha \circ \beta$ by $$\alpha \circ
\beta=\dim \SI(Q,\beta)_{\langle \alpha,\cdot \rangle}=\dim
\SI(Q,\alpha)_{-\langle \cdot, \beta \rangle}.$$

\begin{remark}\label{sat}
As a direct consequence of the saturation property for effective
weights and the above reciprocity property, we have
$$\alpha \circ
\beta \neq 0 \Longleftrightarrow r\alpha \circ s \beta \neq 0,
\forall ~r,~s \geq 1.$$ We also have the following rather trivial
fact
$$\alpha_1 \circ \beta \neq
0,~\alpha_2 \circ \beta \neq 0 \Longrightarrow (\alpha_1 +
\alpha_2) \circ \beta \neq 0.$$ Indeed, multiplying a non-zero
semi-invariant of weight $\langle \alpha_1, \cdot \rangle$ with a
non-zero semi-invariant of weight  $\langle \alpha_2, \cdot
\rangle,$ we obtain a non-zero semi-invariant of weight $ \langle
\alpha_1+\alpha_2, \cdot \rangle.$ Similarly, we have
$$\alpha \circ \beta_1 \neq 0,~\alpha \circ \beta_2 \neq 0
\Longrightarrow \alpha \circ (\beta_1 + \beta_2) \neq 0.$$
\end{remark}

\subsection{$\sigma$-semi-stability}
We have seen that the action of $\GL(\beta)$ on the representation
space $\Rep(Q, \beta)$ gives no interesting quotient varieties. By
twisting the action of $\GL(\beta)$ by means of a weight $\sigma$,
one can obtain plenty of non-trivial semi-invariants. In this way,
King \cite{K} developed a very useful version of GIT to construct
a stability structure for finite dimensional algebras. Let $\sigma
\in \ZZ^{Q_0}$ be a weight such that $\sigma(\beta)=0$.

The following numerical criterion for $\sigma$-(semi-)stability is
due to King \cite{K}. Actually, the original criterion differs
from the one in the theorem below by a sign. This is essentially
because in King's paper \cite{K}, a semi-invariant of weigh
$-\sigma$ is for us a semi-invariant of weight $\sigma$.

\begin{theorem} \label{King}
Suppose that $Q$ is a quiver, $\beta$ a dimension vector,  and $W
\in \Rep(Q,\beta)$. Then:
\begin{enumerate}
\renewcommand{\theenumi}{\arabic{enumi}}
\item $W$ is $\sigma$-semi-stable if and only if for every
subrepresentation $V$ of $W$ we have $$\sigma(\underline d_V) \leq
0;$$ \item $W$ is $\sigma$-stable if and only if for every proper
subrepresentation $V$ of $W$ we have
$$\sigma(\underline d_V) < 0.$$
\end{enumerate}
\end{theorem}

We say that \emph{$\beta$ is $\sigma$-(semi-)stable} if there
exists a $\sigma$-(semi-)stable representation in $\Rep(Q,\beta)$.

\begin{remark}
With this description of $\sigma$-(semi-)stable representations,
one can define the full subcategory of $\Rep(Q)$ consisting of all
$\sigma$-semi-stable representations, including the zero one. Note
that in this abelian category the simple objects are exactly the
$\sigma$-stable representations. Moreover, it can be proved that
this subcategory is Artinian and Noetherian and hence any $\sigma$
-semi-stable representation has a Jordan-Holder filtration with
factors $\sigma$-stable.
\end{remark}

\subsection{General representations of quivers} We will use the
language of general representations of quivers developed by
Schofield \cite{S1} to find necessary and sufficient inequalities
for the non-vanishing of $\dim \SI( Q, \beta)_{\sigma}.$ Let
$\alpha, \beta$ be two dimension vectors. We define the generic
$\ext (\alpha, \beta)$ to be
$$
\ext (\alpha, \beta)= \min \{  \dim \Ext^1_Q(V,W) \mid (V, W) \in
\Rep(Q, \alpha) \times \Rep(Q, \beta) \}.
$$

We write $\alpha \hookrightarrow\beta$ if every representation of
dimension vector $\beta$ has a subrepresentation of dimension
vector $\alpha$. We write $\beta \twoheadrightarrow \alpha $ if
every representation of dimension vector $\beta$ has a quotient
representation of dimension vector $\alpha.$ In other words, we
have that $\alpha \hookrightarrow\beta$ if and only if $\beta
\twoheadrightarrow \beta-\alpha.$ The following lemma follows
immediately from the definition.

\begin{lemma}\label{gensublemma}
Let $\alpha, \beta,\alpha_1, \alpha_2, \beta_1, \beta_2$ be
dimension vectors.
\begin{enumerate}
\renewcommand{\theenumi}{\arabic{enumi}}
\item If $\alpha_2 \hookrightarrow\alpha_1$ and $\alpha_1
\hookrightarrow\alpha$ then $\alpha_2 \hookrightarrow\alpha.$

\item If $\beta \twoheadrightarrow \beta_1$ and $\beta_1
\twoheadrightarrow \beta_2$ then $\beta \twoheadrightarrow
\beta_2.$
\end{enumerate}
\end{lemma}

The next result was proved by Schofield \cite[Theorem 3.3]{S1}.

\begin{lemma}\label{gensub}
Let $\alpha, \beta$ be two dimension vectors. Then the following
are equivalent:
\begin{enumerate}
\renewcommand{\theenumi}{\arabic{enumi}}

\item$\alpha \hookrightarrow\alpha+\beta;$

\item $\ext (\alpha,\beta)=0.$
\end{enumerate}
\end{lemma}

Now, we can give a first description of the set $ \Sigma (Q,
\beta).$

\begin{theorem}\label{semi1}
Let $Q$ be a quiver and $\beta$ be a dimension vector. If $\sigma
= \langle \alpha, \cdot \rangle \in \ZZ^{Q_0}$ is a weight then
the following statements are equivalent:
\begin{enumerate}
\renewcommand{\theenumi}{\arabic{enumi}}
\item $\dim \SI(Q, \beta)_{\sigma} \neq 0,$ i.e., $\sigma \in
\Sigma(Q, \beta);$
\item $\sigma(\beta) = 0$ and $\sigma(\beta')
\leq 0$ for all $\beta'\hookrightarrow \beta;$

\item $\alpha$ must be a dimension vector, $\sigma( \beta) = 0$
and $\alpha \hookrightarrow \alpha+\beta.$
\end{enumerate}
\end{theorem}

\begin{proof}
The equivalence of $(1)$ and $(2)$ is proved in \cite[Theorem
3]{DW1}. It is a direct consequence of the Schofield's \cite{S1}
computation of $\ext (\alpha, \beta)$ and the spanning theorem for
semi-invariants. In the same paper \cite{DW1}, it was noticed that
$\SI (Q, \beta)_{\sigma} \neq 0$ is equivalent to $\alpha$ being a
dimension vector, $\ext (\alpha, \beta)=0$ and $\sigma (\beta) =
\langle \alpha, \beta \rangle = 0$. Hence the equivalence of $(1)$
and $(3)$ follows now from Lemma \ref{gensub}.
\end{proof}

\begin{remark}
It turns out that some of the necessary and sufficient linear
homogeneous inequalities obtained in Theorem \ref{semi1}{(2)} are
redundant. In the next subsection, we will see how one can find a
minimal list of such inequalities.
\end{remark}

A representation $V$ is said to be Schur if $\End_{Q}(V) = \CC.$
We say that a dimension vector $\beta$ is a \emph{Schur} root if
there exists a Schur representation $V$ of dimension vector
$\beta.$ We end this subsection with a description of Schur roots
in terms of $\sigma$- stability.

\begin{theorem}\cite[Theorem 6.1]{S1}\label{stabschur}
Let $Q$ be a quiver and $\beta$ a dimension vector. Then the
following are equivalent:
\begin{enumerate}
\renewcommand{\theenumi}{\arabic{enumi}}
\item $\beta$ is a Schur root; \item $\sigma_{\beta}(\beta')<0,
\forall~ \beta'\hookrightarrow\beta,~\beta'\neq 0, \beta,$ where
$\sigma_{\beta}=\langle \beta,\cdot\rangle-\langle \cdot,\beta
\rangle.$
\end{enumerate}
\end{theorem}

\subsection{$\sigma$-stable decomposition: facets of the cone
$C(Q,\beta)$ of effective weights}Let $\mathbb H(\beta)=\lbrace
\sigma \in \mathbb R^{Q_0} \mid \sigma(\beta)=0 \rbrace$. Consider
the following rational convex polyhedral cone $$ C(Q,\beta)=
\lbrace \sigma \in \mathbb H(\beta) \mid \sigma(\beta') \leq 0
\text{ for all } \beta' \hookrightarrow \beta\rbrace.$$ We call
$C(Q,\beta)$ the \emph{cone of effective weights} associated to
the quiver setting $(Q,\beta)$. Note that $C(Q,\beta)\bigcap
\mathbb Z^{Q_0}=\Sigma(Q,\beta)$ and the dimension of this cone is
at most $N-1$, where $N = |Q_0|$ is the number of vertices of $Q$.

\begin{lemma}\label{interface}
Let $Q$ be a quiver and let $\beta, \gamma_1, \gamma_2, \gamma$ be
dimension vectors.
\begin{enumerate}
\renewcommand{\theenumi}{\arabic{enumi}}
\item Suppose that $\gamma_1+\gamma_2=\beta $ and
$\gamma_1\hookrightarrow \beta$. Then
$$
C(Q,\gamma_1)\bigcap C(Q,\gamma_2)=\mathbb H(\gamma_1)\bigcap
C(Q,\beta).
$$
\item If $c\geq 1$ is a positive integer then
$$
C(Q,c\gamma)=C(Q,\gamma).
$$
\end{enumerate}
\end{lemma}
\begin{proof}
$(1)$ Suppose $\sigma \in C(Q,\gamma_1)\bigcap C(Q,\gamma_2)$ is a
lattice point. Let $\alpha$ be the dimension vector such that
$\sigma=\langle \alpha,\cdot \rangle.$ From Remark \ref{sat}, we
have $$\SI(Q, \alpha)_{-\langle \cdot,\gamma_1\rangle+-\langle
\cdot,\gamma_2\rangle}\neq 0.$$ Using again the reciprocity
property we obtain that $\sigma \in C(Q, \beta)$ and so
$$C(Q,\gamma_1)\bigcap C(Q,\gamma_2)\subseteq \mathbb
H(\gamma_1)\bigcap C(Q,\beta).$$

For the other inclusion, pick a lattice point $\sigma \in \mathbb
H(\gamma_1)\bigcap C(Q,\beta).$ If $\gamma\hookrightarrow\gamma_1$
then from $ \gamma_1 \hookrightarrow \beta$ follows that $\gamma
\hookrightarrow \beta.$ As $\sigma \in \Sigma(Q,\beta),$ we have
that $\sigma(\gamma)\leq 0$ by Theorem \ref{semi1}. This shows
that $\sigma \in C(Q,\gamma_1).$ Now let us assume that
$\delta\hookrightarrow\gamma_2.$ Since $\beta\twoheadrightarrow
\gamma_2$ it follows from transitivity that
$\beta\twoheadrightarrow \gamma_2-\delta$ which is equivalent to
$\gamma_1+\delta\hookrightarrow \beta.$ But this implies
$\sigma(\gamma_1+\delta)=\sigma(\delta)\leq 0.$ We have shown that
$\sigma \in C(Q, \gamma_2),$ as well. Hence
$$\mathbb H(\gamma_1)\bigcap C(Q,\beta) \subseteq
C(Q,\gamma_1)\bigcap C(Q,\gamma_2).$$

$(2)$ This part follows from the reciprocity property (Lemma
\ref{reciprocity}) and the saturation property for effective
weights (Proposition \ref{saturation}).
\end{proof}

An interesting question is to describe the faces of $C(Q,\beta)$.
A useful tool in this direction is the notion of $\sigma$-stable
decomposition for dimension vectors introduced in \cite{DW2}.

Let $\beta$ be a $\sigma$-semi-stable dimension vector. We say
that $$\beta=\beta_1 \pp \beta_2 \pp \ldots \pp \beta_s$$ is the
\emph{$\sigma$-stable decomposition} of $\beta$ if a general
representation in $\Rep(Q,\beta)$ has a Jordan-Holder filtration
(in the full subcategory of $\sigma$-semi-stable representations)
with factors of dimension $\beta_1, \ldots ,\beta_s$ (in some
order). We write $c\cdot\beta$ instead of $\beta \pp \beta \pp
\ldots \pp \beta$ ($c$ times).

\begin{proposition}[\cite{DW2}]\label{basicsigma}
Assume that $\beta$ is a $\sigma$-semi-stable dimension vector. If
$\beta=c_1\cdot \beta_1 \pp c_2\cdot \beta_2 \pp \ldots \pp c_u
\cdot \beta_u$ is the $\sigma$-stable decomposition of $\beta$
with the dimension vectors $\beta_i$ distinct then:

\begin{enumerate}
\renewcommand{\theenumi}{\arabic{enumi}}
\item all $\beta_i$ are Schur roots;

\item if $\langle \beta_i, \beta_i \rangle <0$ then $c_i=1;$

\item after rearranging we can assume that $\beta_i \circ
\beta_j=1$ for all $i<j;$

\item $\beta_1, \dots, \beta_u$ are linearly independent.
\end{enumerate}
\end{proposition}

The relationship between the facets of the cone $C(Q,\beta)$ and
the $\sigma$-stable decomposition is described in the following
Proposition. A stronger form of this Proposition can be found in
\cite[Section 6]{DW2}.

\begin{proposition}\label{descfacets}
Let $Q$ be a quiver with $N$ vertices and let us assume that
$\beta$ is a Schur root.Then
\begin{enumerate}
\renewcommand{\theenumi}{\arabic{enumi}}
\item $\dim C(Q,\beta)=N-1.$

\item $\sigma \in C(Q,\beta)$ if and only if $\sigma(\beta)=0$ and
$\sigma(\beta_1)\leq 0$ for every decomposition $\beta=c_1\beta_1
+ c_2\beta_2$ with $\beta_1,\beta_2$ Schur roots,~
$\beta_1\circ\beta_2=1$ and $c_i=1$ whenever $\langle \beta_i,
\beta_i \rangle <0.$
\end{enumerate}
\end{proposition}

\begin{proof}
$(1)$ Let $C(Q,\beta)^0$ be the open subset of $\mathbb H(\beta)$
defined by
$$
C(Q,\beta)^0=\lbrace \sigma \in \mathbb H(\beta) \mid
\sigma(\beta') < 0 \text{ for all } \beta' \hookrightarrow \beta,
\beta' \neq 0, \beta \rbrace.$$ Since $\beta$ is a Schur root it
follows from Theorem \ref{stabschur} that $\sigma_{\beta} \in
C(Q,\beta)^0$ and hence $C(Q,\beta)^0$ is a non-empty open subset
in $\mathbb H(\beta).$ Consequently, $C(Q,\beta)$ has dimension
$N-1.$

$(2)$ Let $\mathcal F$ be a face of $C(Q,\beta)$ of dimension
$N-2$ and let $\sigma$ be a lattice point in the relative interior
of $\mathcal F.$ Suppose that $$\beta=c_1\cdot \beta_1 \pp
c_2\cdot \beta_2 \pp \dots \pp c_u \cdot \beta_u$$ is the
$\sigma$-stable decomposition of $\beta$ with $\beta_1,
\dots,\beta_u$ as in Proposition \ref{basicsigma}. If
$$\gamma_i=c_1\beta_1 + c_2\beta_2+ \dots + c_i\beta_i$$ for every
$1\leq i \leq u$ then it follows from Remark \ref{sat} and Theorem
\ref{semi1} that $\gamma_{i}\hookrightarrow \beta$ for every
$1\leq i\leq u.$ This shows that the $-\gamma_i$'s, viewed as
linear forms on $\mathbb R^{Q_0}$, are in the dual cone of
$C(Q,\beta)$ and hence
$$\mathbb H(-\gamma_1)\bigcap \dots \bigcap \mathbb
H(-\gamma_u)\bigcap C(Q,\beta)$$ is a face, denoted by $\mathcal
F',$ of $C(Q,\beta).$ As $\mathcal F'$ is a face of $C(Q,\beta)$
containing a relative interior point $\sigma$ of another face
$\mathcal F$ it follows that $\mathcal F \subseteq \mathcal F'.$
We have that $\gamma_1, \dots, \gamma_u$ are linearly independent
as $\beta_1, \dots, \beta_u$ have this property. Thus, the dual
face (in $(\mathbb R^{Q_0})^{*}$) of $\mathcal F'$ has dimension
at least $u$ and so $u\leq 2.$

If $u=1$ then $C(Q,\beta)=C(Q,\beta_1)$ by Lemma \ref{interface}.
As $\beta_1$ is $\sigma$-stable and using Theorem \ref{King} we
obtain that $\sigma$ must lie in the relative interior of
$C(Q,\beta).$ But this is a contradiction with the fact that
$\sigma$ lies on a proper face $\mathcal F$ of $C(Q, \beta).$

Therefore, our facet $\mathcal F$ has to be of the form
$$
\mathbb H(-\beta_1)\bigcap C(Q,\beta) = C(Q,\beta_1)\bigcap
C(Q,\beta_2)
$$
for some Schur roots $\beta_1,\beta_2$ with
$\beta_1\circ\beta_2=1$ and $\beta=c_1\beta_1+c_2\beta_2$ with
$c_1, c_2$ as in Proposition \ref{basicsigma}. This description of
the facets of the cone $C(Q,\beta)$ together with Theorem
\ref{semi1} clearly imply $(2).$
\end{proof}

\begin{remark} \label{longlist}
Note that in Proposition \ref{descfacets}{(2)}, we can replace
$\beta_1\circ\beta_2 = 1$ with $\beta_1\circ\beta_2 \neq 0.$ Of
course, in this case we get a longer list of necessary and
sufficient inequalities.
\end{remark}

\begin{remark} \label{facetsrmk}
In \cite{DW2} (see also \cite{CC}), it has been conjectured that
the list of linear homogeneous inequalities from Proposition
\ref{descfacets}{(2)} is minimal, i.e., the facets of $C(Q,\beta)$
when $\beta$ is a Schur root are in one-to-one correspondence with
the set of pairs $(\beta_1, \beta_2)$ where $\beta_1$ and
$\beta_2$ are as in Proposition \ref{descfacets}{(2)}.
\end{remark}

\section{The facets of the cone associated to the generalized flag
quiver}\label{facetsgenflag} We use those methods from Section
\ref{coneabstract} to describe the facets of the cone of effective
weights associated to the generalized flag quiver setting.

Throughout this section the only quiver we will be dealing with is
the generalized flag quiver setting $(Q, \beta )$ from Section
\ref{flagquiver}. For the convenience of the reader, we briefly
recall this set up. The quiver $Q$ has $m-2$ central vertices with
$m$ equioriented $\mathbb A_n$ quivers (or flags) $\mathcal
F(1),\dots,\mathcal F(m)$ attached to them. The dimension vector
$\beta$ is defined by $\beta(j, i) = j$ for all $j \in \{1, \dots,
n\}$ and $i \in \{1, \dots, m\}.$

First, let us prove a simple Lemma.

\begin{lemma}
The dimension vector $\beta$ is a Schur root.
\end{lemma}

\begin{proof}
Note that the dimension vector $\beta$ is indivisible, meaning
that the greatest common divisor of its coordinates is one. Next,
let us assume that either  $n=2, ~m \geq 4$ or $n \geq 3.$ If this
is the case then $\beta$ lies in the so called \emph{fundamental
region,} i.e., the support of $\beta$ is a connected graph and
$\langle \varepsilon_i,\beta \rangle + \langle \beta,
\varepsilon_i \rangle \leq 0,$ for all vertices $i \in Q_0.$ It
follows now from a result of Kac \cite[Theorem B(d)]{Kac} that
$\beta$ is a Schur root. If either $n=2,~ m=3$ or $n=1$ then
$\beta$ is actually a real Schur root.
\end{proof}

Now, let $\mathcal D$ be the set of all dimension vectors
$\beta_1$ that take one of the following forms:
\begin{enumerate}
\renewcommand{\theenumi}{\arabic{enumi}}
\item $\beta_1 = \varepsilon_{(j, 2i+1)}$ or $\beta_1 = \beta -
\varepsilon_{(j, 2i)},$ for $1 \leq j \leq n-1$ (call such a
dimension vector \emph{trivial});
\newline
or
\item $\beta_1$ is weakly increasing with jumps of at most one
along the $m$ flags, $\beta_1 \neq \beta$ and $\beta_1 \circ
(\beta -\beta_1) = 1.$
\end{enumerate}
Note that if $\beta_1$ is in $\mathcal D$ then $\beta_1
\hookrightarrow \beta$ and hence $-\beta_1$ is in the dual of the
cone $C(Q, \beta).$

\begin{lemma}\label{mainlemma} Keep notation as above. If $\mathcal F$ is a facet
of $C(Q, \beta)$ then it has to be of the form
$$
\mathcal F=\mathbb H(\beta_1)\bigcap C(Q, \beta),
$$
for some $\beta_1$ in $\mathcal D.$
\end{lemma}

\begin{proof}
From Proposition \ref{descfacets} it follows that there are two
Schur roots $\beta_1$ and $\beta_2$ such that $$\mathcal F=\mathbb
H(\beta_1)\bigcap C(Q, \beta)$$ with $\beta_1\circ\beta_2=1$ and
$\beta=c_1\beta_1+c_2\beta_2$ for some $c_1,c_2 \geq 1.$

Now let us assume that $\beta_1$ is not trivial. In this case, we
will show that $\beta_1$ is weakly increasing with jumps of at
most one along the flags. Let us denote $c_1\beta_1=\beta',
c_2\beta_2=\beta''$. Since $\beta'\circ\beta'' \neq 0$ it follows
from Theorem \ref{semi1} that any representation of dimension
vector $\beta$ has a subrepresentation of dimension vector
$\beta'.$ Therefore, $\beta'$ must be weakly increasing along each
flag going in and it has jumps of at most one along each flag
going out.

Next, we will show that $\beta'$ has jumps of at most one along
each flag $\mathcal F(i)$ going in a central vertex and $\beta'$
is weakly increasing along each flag $\mathcal F(i)$ going out
from a central vertex. For simplicity, let us write
\begin{diagram} \mathcal F(i):
1&\rTo&2&\cdots&n-1&\rTo&n,
\end{diagram}
for a flag going in its central vertex $(n,i)$ (i.e. $i$ is even).
Assume to the contrary that there is an $l\in \{1, \dots, n-1\}$
such that $\beta'(l+1)>\beta'(l)+1.$ Then $\beta ''(l+1) < \beta
''(l)$ which implies that $\varepsilon_l\hookrightarrow\beta ''.$
Since $\beta''$ is $\langle\beta',\cdot \rangle$-semi-stable it
follows that $\langle\beta',\epsilon_l\rangle\leq 0.$ So,
$\beta'(l)\leq \beta'(l-1)$ and hence $\beta'(l)=\beta'(l-1)$ or
$\beta''(l)=\beta''(l-1)+1.$ This shows that $c_2=1$ and $\beta''
- \varepsilon_l\hookrightarrow\beta ''.$ From the fact that $\beta
''(=\beta_2)$ is a Schur root and Theorem \ref{stabschur} we
obtain that $\beta ''$ is $\sigma_{\beta ''}$-stable. Since
$\varepsilon_l\hookrightarrow\beta
''$,$\beta"-\varepsilon_l\hookrightarrow\beta ''$ and $\beta ''
\neq \epsilon_l$ it follows $\langle\beta
'',\epsilon_l\rangle-\langle\beta '',\varepsilon_l\rangle<0$ and
$\langle \beta '',\beta ''-\varepsilon_l\rangle-\langle \beta '' -
\varepsilon_l, \beta '' \rangle < 0.$ But this is a contradiction.
We have just proved that $\beta'$ has jumps of at most one along
each flag going in. Similarly, one can show that $\beta'$ has to
be weakly increasing along each flag going out.

Now, let us show that $c_1=c_2=1$. Since $\beta'=c_1\beta_1$ has
jumps of at most one along each flag, we obtain $0\leq
c_1(\beta_1(l+1, i)-\beta_1(l, i))\leq 1$ for all $l\in \{1,
\dots, n-1\}$ and $i \in \{1, \dots, m\}.$ If there are $l,i$ such
that $\beta_1(l+1, i)-\beta_1(l, i)\neq 0$ then $c_1=1$.
Otherwise, there is an $i$ such that $\beta'(1, i)=1$ and so
$c_1=1$. Similarly, one can show $c_2=1$.

In conclusion, $\beta=\beta_1+\beta_2$ with $\beta_1$ weakly
increasing with jumps of at most one along the $m$ flags. So,
$\beta_1 \in \mathcal D$ and this finishes the proof.
\end{proof}

\begin{lemma} \label{maincoro}
Let $\sigma \in \mathbb H( \beta).$ Then
$$\sigma \in C(Q, \beta)$$ if and only if the following
are true
\begin{enumerate}
\renewcommand{\theenumi}{\arabic{enumi}}
\item (\emph{chamber inequalities}) $(-1)^i\sigma(\varepsilon_{(j,
i)}) \geq 0, ~\forall ~ 1 \leq j \leq n-1,~\forall ~1 \leq i \leq
m.$

\item (\emph{regular inequalities}) $\sigma(\beta_1)\leq 0$ for
every $\beta_1 \neq \beta$ weakly increasing with jumps of at most
one along the $m$ flags and $\beta_1 \circ (\beta - \beta_1) = 1.$
\end{enumerate}
\end{lemma}

\begin{proof}
Let us assume that $\sigma \in \mathbb H(\beta)$ satisfies the
chamber and regular inequalities. Then the description of the
facets of $C(Q, \beta)$ given in Lemma \ref{mainlemma} shows that
$\sigma \in C(Q, \beta).$

Conversely, let $\sigma \in C(Q, \beta).$ We clearly have $\sigma(
\beta_1) \leq 0$ for every $\beta_1 \in \mathcal D$ by Theorem
\ref{semi1}. But this is equivalent to $(1)$ and $(2).$
\end{proof}

\begin{remark}
Let $\sigma_{\lambda}$ be the weight defined by the equations
$(\ref{thewt1})-(\ref{thewt4})$ in Section \ref{flagquiver}. Then
by definition we have that
$$
\sigma_{\lambda}(\varepsilon_{(j,
i)})=(-1)^i(\lambda_j(i)-\lambda_{j+1}(i)), \forall ~1 \leq j \leq
n-1, ~\forall ~1 \leq i \leq m.
$$
Consequently, the \emph{chamber inequalities} just tell us that
the $\lambda(i)$ are weakly decreasing sequences. This is
something that we will always assume.
\end{remark}
%%%%%%%%%%%%%%%%%%%%%%%%%%%%%%%%%%%%%%%%%%%%%%%%%%%%%%%%%%%%%%%%%%%%%%%%%%%%
\begin{example}\label{ex}
For $m=4$ and $n=2$, there are exactly $9$ dimension vectors
$\beta_1$ that satisfy the second condition in Lemma
\ref{maincoro}. It turns out that exactly one of the $9$ pairs
gives us a redundant inequality. Next we find the necessary and
sufficient inequalities for $\sigma_{\lambda}$ to be in $C( Q,
\beta).$

For $$\beta_1 = \beta =
\begin{matrix}
 & 2 & 2 \\
1 & 1 & 1 & 1
\end{matrix},$$ we must have the identity $\sigma_{\lambda}(\beta) =
0,$ i.e.,
$$
|\lambda(1)|+|\lambda(3)|=|\lambda(2)|+|\lambda(4)|.
$$
%%%%%%%%%%%%%%%%%%%%%%%%%%%%%%%%%%%%%%%%%%%%%%%%%%%%%%%%%%
For $\beta_1=
\begin{matrix}
 & 1 & 2 \\
0 & 0 & 1 & 1
\end{matrix}$
and
$ \beta_1=
\begin{matrix}
 & 1 & 2 \\
1 & 1 & 1 & 1
\end{matrix},$ we have the inequalities
$$
\lambda_2(2)+|\lambda(4)|\leq \lambda_2(1)+|\lambda(3)|,
$$
and
$$
\lambda_1(2)+|\lambda(4)|\leq \lambda_1(1)+|\lambda(3)|.
$$
%%%%%%%%%%%%%%%%%%%%%%%%%%%%%%%%%%%%%%%%%%%%%%%%%%%%%%%%%%%
For $\beta_1=
\begin{matrix}
 & 0 & 1 \\
0 & 0 & 1 & 1
\end{matrix}
$
and
$\beta_1=
\begin{matrix}
 & 0 & 1 \\
0 & 0 & 0 & 0
\end{matrix},$ we have the inequalities
$$
\lambda_1(4)\leq \lambda_1(3),\text{~and~} \lambda_2(4)\leq
\lambda_2(3).
$$
%%%%%%%%%%%%%%%%%%%%%%%%%%%%%%%%%%%%%%%%%%%%%%%%%%%%%%%%%%%%
For $\beta_1=
\begin{matrix}
 & 1 & 1 \\
1 & 0 & 1 & 1
\end{matrix}
$  and  $\beta_1=
\begin{matrix}
 & 1 & 1 \\
1 & 1 & 1 & 0
\end{matrix},$ we have the inequalities
$$
\lambda_2(2)+\lambda_1(4)\leq \lambda_1(1)+\lambda_1(3),
$$
and
$$
\lambda_1(2)+\lambda_2(4)\leq \lambda_1(1)+\lambda_1(3).
$$
%%%%%%%%%%%%%%%%%%%%%%%%%%%%%%%%%%%%%%%%%%%%%%%%%%%%%%%%%%%%%%%%
For $\beta_1=
\begin{matrix}
 & 1 & 1 \\
0 & 0 & 1 & 0
\end{matrix}
$ and $\beta_1=
\begin{matrix}
 & 1 & 1 \\
1 & 0 & 0 & 0
\end{matrix},$ we have the inequalities
$$
\lambda_2(2)+\lambda_2(4)\leq \lambda_2(1)+\lambda_1(3),
$$
and
$$
\lambda_2(2)+\lambda_2(4)\leq \lambda_1(1)+\lambda_2(3).
$$\\
%%%%%%%%%%%%%%%%%%%%%%%%%%%%%%%%%%%%%%%%%%%%%%%%%%%%%%%%%%%%%5
For $$\beta_1=
\begin{matrix}
 & 0 & 2 \\
0 & 0 & 1 & 1
\end{matrix},$$ we obtain the only redundant inequality
$$
|\lambda(4)| \leq |\lambda(3)|.
$$
\end{example}

%%%%%%%%%%%%%%%%%%%%%%%%%%%%%%%%%%%%%%%%%%%%%%%%%%%%%%%%%%%%%%%

\section{The Horn type inequalities}\label{proofs}
Our goal in this section is to give a closed form to the
polyhedral inequalities that we obtained in Lemma \ref{maincoro}.
The quiver setting that we work with is again the generalized flag
quiver from Section \ref{flagquiver}.

First, let us describe the dimension vectors $\beta_1$ that define
the regular inequalities from Lemma \ref{maincoro}{(2)}. Let
$\beta_1$ be a dimension vector that is weakly increasing with
jumps of at most one along the $m$ flags. We define the following
jump sets
$$
I_i = \{l \mid \beta_1(l, i)>\beta_1(l-1, i),1 \leq l\leq n \},
$$
with the convention that $\beta_1(0, i)=0$ for all $i \in \{1,
\dots, m \}.$ We also denote $\beta_1$ by $\beta_I.$

Note also that $|I_i|=\beta_I(n, i)$ for all $i\in\{1, \dots, m
\}.$ Therefore, $|I_1| = |I_2| = \beta_I(2)$ and $|I_{m-1}| =
|I_m| = \beta_I(m-1).$

Conversely, it is clear that each $m$-tuple $I=(I_1, \dots, I_m)$
of subsets of the set $\{1, \dots, n \}$ with  $|I_1| = |I_2|,$
$|I_{m-1}| = |I_m|,$ uniquely determines the dimension vector
$\beta_I.$ Indeed, if
$$
I_i = \{z_1(i)< \dots < z_r(i)\},
$$
we have that
$$
\beta_I(k, i)=j-1, \forall~ z_{j-1}(i)\leq  k < z_j(i), \forall~ 1
\leq j \leq r+1,
$$
with the convention that $z_0(i)=0$ and $z_{r+1}(i) = n+1$ for all
$1 \leq i \leq m.$

\begin{definition} \label{theset}
We define $\mathcal S(n, m)$ to be the set consisting of all
$m$-tuples $I=(I_1, \dots, I_m)$ such that $|I_1| = |I_2|,$
$|I_{m-1}| = |I_m|,$ $\beta_I \neq \beta $ and
$$
\beta_I \circ (\beta-\beta_I) =1.
$$
\end{definition}

A further description of the set $\mathcal S(n, m)$ will be given
in Lemma \ref{descset} and Lemma \ref{mainlemma2}.

\begin{proposition}\label{mainprop}
Let $\lambda(1), \dots, \lambda(m)$ be weakly decreasing sequences
of $n$ reals. Then the following are equivalent:
\begin{enumerate}
\renewcommand{\theenumi}{\arabic{enumi}}
\item $\sigma_{\lambda} \in C(Q, \beta);$

\item
$$\sum_{i~even}|\lambda(i)|=\sum_{i
~odd}|\lambda(i)|$$ and
$$\sum_{i~even} \left(\sum_{j \in
I_i}\lambda_j(i)\right) \leq \sum_{i~odd} \left(\sum_{j\in
I_i}\lambda_j(i)\right)$$ for every $m$-tuple $(I_1, \dots, I_m)
\in \mathcal S(n, m).$
\end{enumerate}
\end{proposition}

\begin{proof}
We have seen that the set of all $\beta_1$ occurring in Lemma
\ref{maincoro}{(2)} are exactly those of the form $\beta_I$ with
$I=(I_1, \dots, I_m) \in \mathcal S(n, m).$ Furthermore it is easy
to see that
$$\sigma_{\lambda}(\beta_I)=\sum_{i~even } \left(\sum_{j \in I_i}\lambda_j(i)\right)-\sum_{i~
odd} \left(\sum_{j\in I_i}\lambda_j(i)\right)$$ and
$$
\sigma_{\lambda}(\beta)=\sum_{i~even}|\lambda(i)|-\sum_{i
~odd}|\lambda(i)|.
$$ The Proposition is now an immediate consequence of Lemma \ref{maincoro}.
\end{proof}

\begin{example}
In this example we will work out the case when $n=1.$ Let $d_1,
\dots, d_m$ be $m\geq 3$ positive integers. Then the following are
equivalent:
\begin{enumerate}
\renewcommand{\theenumi}{\arabic{enumi}}
\item There exists a long exact sequence of the form
$$
0 \to (\ZZ/p)^{d_1} \to \cdots \to (\ZZ/p)^{d_m} \to 0.
$$
\item There exists a long exact sequence of the form
$$
0 \to \ZZ/p^{d_1} \to \cdots \to \ZZ/p^{d_m} \to 0.
$$
\item(\textit{Horn type inequalities})
$$
\sum_{j~even} d_j=\sum_{j~odd} d_j
$$
and if $m >3$
$$
\sum_{j~ even,~1\leq j \leq i} d_j \leq \sum_{j~ odd,~1\leq j \leq
i} d_j
$$
and
$$
\sum_{j~ even,~i \leq j \leq m} d_j \leq \sum_{j~ odd,~i\leq j
\leq m} d_j,
$$
for every $i$ odd with $2 \leq i \leq m-2,$ together with $d_m
\leq d_{m-1}$ if $m$ is even.
\end{enumerate}
\end{example}

Indeed, let $\lambda(i)=(d_i), \forall~ 1\leq i \leq m$. The
equivalence of $(1)$ and $(2)$ follows from
$$
f(\lambda(1), \dots, \lambda(m))\neq 0 \Longleftrightarrow
f(\lambda'(1), \dots, \lambda'(m))\neq 0.
$$

To prove the equivalence $(2) \Longleftrightarrow (3),$ we
explicitly describe the facets of the cone $C(Q, \beta),$ where
$Q$ is the generalized quiver when $n=1.$ When $m=3$, the only
inequality is $d_2=d_1+d_3$. Let us assume that $m\geq 4$. In this
case, our quiver $Q$ is an alternating type $\mathbb A_{m-2}$
quiver with $m-2$ vertices such that $2$ is a source, $3$ is a
sink and so on. For example if $m$ is odd then our generalized
flag quiver becomes

\begin{diagram}
 2& \rTo & 3 & \lTo & \cdots & m-2 &\lTo & m-1,
\end{diagram}

First, let $\beta_1, \beta_2$ be two Schur roots (i.e. positive
roots of type $\mathbb A$) such that
$\beta_1+\beta_2=\beta=(1,\dots, 1)$ and $\langle
\beta_1,\beta_2\rangle =0$. Then it is easy to see that
$$
\beta_1=(1,\dots,1 ,0,\dots, 0)~ \text{or}~ \beta_1=(0,\dots, 0,1,
\dots,1),
$$
with $\supp(\beta_1)=\{2, \dots, i\}~ \text{or}~ \{i, \dots,
m-1\}$ and $2\leq i \leq m-1$ odd.  To find a minimal list of
necessary and sufficient inequalities, we will focus on those
$m$-tuples $I=(I_1, \dots, I_m)\in \mathcal S$ for which the
corresponding dimension vectors $\beta_I,~\beta-\beta_I$ are Schur
roots. If this the case, we must have that
$$I_j =
\begin{cases}
\{1\} & \text{if $1\leq j \leq i$} \\
\emptyset & \text{if $i<j \leq m$ }
\end{cases}
$$
or
$$I_j =
\begin{cases}
\emptyset & \text{if $1\leq j<i$} \\
\{1\} & \text{if $i\leq j \leq m$, }
\end{cases}
$$
where $2 \leq i \leq m-2$ is odd. If $m$ is even, there is one
more possibility, namely $\beta_1 = (0, \dots, 0, 1).$ In this
case, $I_1 = \dots = I_{m-2}=\emptyset$ and $I_{m-1} = I{m} =
\{1\}.$ For all such tuples $I,$ we also have that $\beta_I \circ
(\beta - \beta_I) = 1.$ This way, we obtain the equivalence of
$(2)$ and $(3).$ Note that the list of inequalities obtained is
minimal.
\newline

Now, let us show that $\mathcal S(n, m)$ can be described in terms
of the generalized Littlewood-Richardson coefficients. For
convenience, let us recall some of the notation from Section
\ref{intro}. Let $(I_1, \dots, I_m)$ be an $m$-tuple of subsets of
$\{1, \dots, n \}$ such that at least one of them has cardinality
at most $n-1.$ We define the following weakly decreasing sequences
of integers (using conjugate partitions):

$\underline\lambda(I_1)=\lambda'(I_1),~~~~\underline\lambda(I_m)=
\begin{cases}
\lambda'(I_m) & \text{if $m$ is odd} \\
\lambda'(I_m\setminus \{n\}) & \text{if $m$ is even},
\end{cases}
$

and for $2\leq i\leq m-1$ \newline

$\underline\lambda(I_i)=
\begin{cases}
\lambda'(I_i) & \text{if $i$ is even} \\
\lambda'(I_i)-((|I_i|-|I_{i+1}|-|I_{i-1}|)^{n-|I_{i}|}) & \text{if
$i \leq m-2$ is odd}\\
\lambda'(I_i)-((|I_{m-1}|-|I_{m-2}|-|I_{m}\setminus
\{n\}|)^{n-|I_{i}|}) & \text{if $i = m-1$ is odd.}
\end{cases}
$

\begin{lemma} \label{descset}
The set $\mathcal S(n, m)$ \index{$\mathcal S(n, m)$} consists of
all $m$-tuples $I = (I_1, \dots, I_m)$ such that:
\begin{enumerate}
\renewcommand{\theenumi}{\alph{enumi}}

\item $|I_1| = |I_2|;$

\item $|I_{m-1}| = |I_m|;$

\item at least one of the subsets $I_1, \dots, I_m$ has
cardinality $< n;$

\item $\underline{\lambda}(I_i)$ is a partition, $\forall~1 \leq i
\leq m;$

\item we have
$$f(\underline{\lambda}(I_1), \dots, \underline{\lambda}(I_m)) =
1.$$
\end{enumerate}
\end{lemma}

\begin{proof}
Let $I=(I_1, \dots, I_m)$ be an $m$-tuple in $\mathcal S(n, m).$
By definition, we know that $(a)$ and $(b)$ are satisfied.

Let us denote $\beta_I=\beta_1$ and $\beta- \beta_I=\beta_2.$

$(c)$ If $\min_{1 \leq i \leq m}|I_i|=n$ then we would have
$\beta_1=\beta$ which is not allowed.

$(d),~(e)$ We compute the dimension
$\beta_1\circ\beta_2=\dim\SI(Q,\beta_2)_{\langle \beta_1,\cdot
\rangle}$ using the same arguments as in Lemma \ref{compute} with
$\beta$ replaced by $\beta_2$ and $\sigma$ by $\sigma_1 = \langle
\beta_1,\cdot \rangle.$ Since $\beta_1$ is weakly increasing and
has jumps of at most one along the flags it is easy to see that
$$\sigma_1(l, i)=
\begin{cases}
1 & \text{if $l\in I_i$} \\
0 & \text{otherwise}
\end{cases},$$ \\ for all $l\in\{1, \dots, n-1\}$ and $i$ even and\\
$$\sigma_1(l, i)=
\begin{cases}
-1 & \text{if $l+1\in I_i$} \\
0 & \text{otherwise}
\end{cases},$$ \\for all $l\in\{1, \dots, n-1\}$ and $i$ odd. At the
central vertices $2, \dots, m-1$ the values of $\sigma_1$ are
$$\sigma_1(i)=
\begin{cases}
0 & \text{if $i$ is even and $n \notin I_i$}\\
1 & \text{if $i$ is even and $n \in I_i$}\\
|I_i|-|I_{i+1}|-|I_{i-1}| & \text{if $i \leq m-2$ is odd} \\
|I_{m-1}|-|I_{m-2}|-|I_m\setminus \{n\}| & \text{if $i=m-1$ is
odd.}
\end{cases}
$$
Arguing as in Lemma \ref{compute}, we obtain that
$$
\begin{aligned}
\gamma(1) &
=(\beta_2(n-1, 1)^{-\sigma_1(n-1, 1)}, \dots, \beta_2(1, 1)^{-\sigma_1(1, 1)})', \\
\gamma(m) & =(\beta_2(n-1, m)^{(-1)^m \cdot
\sigma_1(n-1, m)},\dots,\beta_2(1, m)^{(-1)^m \cdot \sigma_1(1, m)})', \\
\gamma(i) & = (\beta_2(n-1, i)^{(-1)^i \cdot \sigma_1(n-1, i)},
\dots, \beta_2(1, i)^{(-1)^i \cdot \sigma_1(1, i)})' + (((-1)^i
\cdot \sigma_1(n, i))^{\beta_2(n, i)}),
\end{aligned}
$$
must be partitions for all $2 \leq i \leq m-1$ and
$$\dim
\SI(Q,\beta_2)_{\sigma_1}=f(\gamma(1), \dots, \gamma(m)).$$
Furthermore, if $I_i = \{z_1(i)< \dots< z_r(i)\}$ then we have
$$\beta_2(z_j(i), i)=z_j(i)-j=\beta_2(z_j(i)-1, i)$$  for all
$j\in\{1, \dots, r\}.$

Therefore, $\gamma(i)=\underline{\lambda}(I_i), 1 \leq i \leq m$
and so $$ f(\underline{\lambda}(I_1), \dots,
\underline{\lambda}(I_m))=1.
$$
We have just proved that if $(I_1, \dots, I_m)$ is in $\mathcal
S(n, m)$ then $(a)-(e)$ are fulfilled.

Conversely, let $I=(I_1, \dots, I_m)$ be an $m$-tuple of subsets
of $\{1, \dots, n \}$ satisfying $(a)-(e).$ Then we can define
$\beta_I$ such that $\beta_I \neq \beta$ and
$$
\beta_I \circ (\beta-\beta_I)=f(\underline{\lambda}(I_1), \dots,
\underline{\lambda}(I_m)) = 1.
$$
Thus, $I=(I_1, \dots, I_m) \in \mathcal S(n, m)$ and so we are
done.
\end{proof}

\begin{proposition}\label{mainpropnou}
Let $\lambda(i)=(\lambda_1(i), \dots, \lambda_n(i)), ~i\in \{1,
\dots, m \}$ be $m$ weakly decreasing sequences of $n$ reals. Then
the following are equivalent:
\begin{enumerate}
\renewcommand{\theenumi}{\arabic{enumi}}
\item $\sigma_{\lambda} \in C(Q, \beta);$

\item the numbers $\lambda_j(i)$ satisfy
$$
\sum_{i~even}|\lambda(i)|=\sum_{i~odd}|\lambda(i)|
$$
together with
$$
(*) \hspace{15pt}\sum_{i~even} \left(\sum_{j \in
I_i}\lambda_j(i)\right) \leq \sum_{i~odd} \left(\sum_{j\in
I_i}\lambda_j(i)\right)
$$
for every $m$-tuple $(I_1, \dots, I_m)$ for which $|I_1| = |I_2|,$
$|I_{m-1}| = |I_m|,$ $\underline{\lambda}(I_i),~ 1 \leq i \leq m$
are partitions and
$$f(\underline{\lambda}(I_1), \dots,
\underline{\lambda}(I_m))\neq 0;$$

\item the numbers $\lambda_j(i)$ satisfy $$
\sum_{i~even}|\lambda(i)|=\sum_{i~odd}|\lambda(i)|
$$ and $(*)$ for
every $m$-tuple $(I_1,\dots, I_m)$ for which $|I_1| = |I_2|,$
$|I_{m-1}| = |I_m|,$ $\underline{\lambda}(I_i),~1 \leq i \leq m$
are partitions and
$$f(\underline{\lambda}(I_1), \dots,
\underline{\lambda}(I_m)) = 1.$$
\end{enumerate}
\end{proposition}

\begin{proof}
The proof follows from Proposition \ref{mainprop}, Lemma
\ref{descset} and Remark \ref{longlist}.
\end{proof}

We end this section with some further remarks on the set $\mathcal
S(n, m).$ The next Lemma gives us constraints on the possible
$m$-tuples $I=(I_1, \dots, I_m)$ of the set $\mathcal S(n, m).$

\begin{lemma}\label{mainlemma2} Let $I=(I_1, \dots, I_m)$ be in $\mathcal S(n, m).$
Then the subsets $I_1, \dots, I_m$ satisfy:
\begin{enumerate}
\renewcommand{\theenumi}{\alph{enumi}}
\item (if $m > 3$) for each $i$ odd, $2\leq i\leq m-2$
$$\max \{|I_{i-1}|, |I_{i+1}|\}\leq |I_i|\leq
|I_{i-1}|+|I_{i+1}|+s_i,$$ where $s_i$ is the smallest $k \in \{0,
\dots, |I_i|\}$ such that $n-k \notin |I_i|;$

\item if $i=m-1$ is odd we have $|I_{m-2}| \leq |I_{m-1}|$ and if
$n\in I_m$ then either $n\in I_{m-1}$ or $I_{m-2}\neq \emptyset.$
\end{enumerate}
\end{lemma}

\begin{proof}
$(a)$ Let us denote $\beta_I=\beta_1$ and $\beta-\beta_I=\beta_2.$
Since $\beta_1\circ\beta_2\neq 0$ we have from Theorem \ref{semi1}
that any representation $V$ of dimension vector
$\beta=\beta_1+\beta_2$ has a subrepresentation of dimension
vector $\beta_1.$ Choose $V$ such that $V(a)$ is invertible for
every main arrow $a$. Then for each $i$ odd, $2\leq i\leq m-1,$ we
clearly have $$\max \{|I_{i-1}|, |I_{i+1}|\}\leq |I_i|.$$

Let us denote $\langle\beta_1,\cdot\rangle$ by $\sigma_1.$ A
necessary condition for $\dim\SI(Q,\beta_2)_{\langle \beta_1,\cdot
\rangle}$ not to be zero is that $\underline\lambda(I_i),$
$\forall ~1 \leq i \leq m$ be partitions, i.e. they must have
non-negative parts.

Suppose that $2\leq i\leq m-1$ is odd and let $s_i$ be the
smallest $k \in \{0, \dots, |I_i|\}$ such that $n-k \notin I_i.$
Then the smallest part of $\lambda'(I_i)$ is exactly $s_i.$

For $2\leq i \leq m-2$ odd, we have seen that
$\underline\lambda(I_i)=\lambda'(I_i)-(\sigma_1(i)^{n-|I_i|}).$ On
the other hand, we know that
$\sigma_1(i)=|I_i|-|I_{i-1}|-|I_{i+1}|$ and the smallest part of
$\lambda'(I_i)$ is precisely $s_i.$ Thus, $\underline\lambda(I_i)
$ is a partition if and only if
$$0 \leq |I_{i-1}|+|I_{i+1}|-|I_i|+s_i.$$

$(b)$ If $i=m-1$ is odd and $n \notin I_m$ then
$$\sigma_1(m-1)=|I_{m-1}|-|I_m|-|I_{m-2}|=-|I_{m-2}|\leq 0$$ in
which case $\underline\lambda(I_{m-1})$ is clearly a partition.

Now let assume that $i=m-1$ is odd and $n \in I_m.$ Then
$$\sigma_1(m-1)=|I_{m-1}|-|I_m|+1-|I_{m-2}|=1-|I_{m-2}|$$ and hence
$\underline\lambda(I_{m-1})$ is a partition when
$$s_{m-1}+|I_{m-2}| \geq 1.$$ So, in this case we must have that
either $n\in I_{m-1}$ or $I_{m-2}\neq \emptyset$.
\end{proof}

\begin{remark} \label{Klycone}
When $m=3,$ the set $\mathcal S(n, 3)$ is just the set of all
triples $(I_1, I_2, I_3)$ of subsets of $\{1, \dots, n\}$ of the
same cardinality $r$ with $r<n$ and
$c_{\lambda(I_1),\lambda(I_3)}^{\lambda(I_2)}=1.$ So, $\mathcal
K(n, 3)$ is indeed the Klyachko's cone. Therefore, in this case we
recover the Horn type inequalities that solve the non-vanishing of
the Littlewood-Richardson coefficients problem and Horn's
conjecture.
\end{remark}

%%%%%%%%%%%%%%%%%%%%%%%%%%%%%%%%%%%%%%%%%%%%%%%%%%%%%%%%%%%%%%%%%%
%%%%%%%%%%%%%%%%%%%%%%%%%%%%%%%%%%%%%%%%%%%%%%%%%%%%%%%%%%%%%%%%%%

\section{Proof of Theorem \ref{mainthm} and Proposition \ref{genKlycone}}\label{proofsnou}

Before we prove our main theorem, we briefly recall the following
moment map description of the cone of effective weights.

\begin{proposition} \cite[Proposition 1.3]{CC3} \label{moment}
Let $Q$ be a quiver without oriented cycles, $\beta$ be a
dimension vector and $\sigma \in \mathbb R^{Q_0}.$ Then the
following statements are equivalent:
\begin{enumerate}
\renewcommand{\theenumi}{\arabic{enumi}}
\item $\sigma \in C(Q,\beta);$

\item  there exists $W=\{W(a)\}_{a \in Q_1} \in \Rep(Q,\beta)$
satisfying
$$(\dag) \hspace{10pt}
\sum_{\buildrel {a \in Q_1} \over {ta=x} }W(a)^*W(a)-\sum_{
\buildrel {a \in Q_1} \over {ha=x}
}W(a)W(a)^*=\sigma(x)\text{Id}_{\beta(x)},
$$
for all $x \in Q_0$, where $W(a)^*$ is the adjoint of $W(a)$ with
respect to the standard Hermitian inner product on $\CC^n.$
\end{enumerate}
\end{proposition}

In what follows, we work with the generalized flag quiver setting
from Section \ref{flagquiver}. To apply Proposition \ref{moment},
we need the following simple linear algebra Lemma:

\begin{lemma} \label{linalg}
Let $\sigma(1), \dots, \sigma(n-1)$ be non-positive real numbers.
Then the following are equivalent:
\begin{enumerate}
\renewcommand{\theenumi}{\arabic{enumi}}
\item there exist $W_i\in \Mat_{i\times (i+1)}(\CC), 1\leq i\leq
n-1$ such that
$$
\begin{aligned}
W_i\cdot W_i^{\star} - W_{i-1}^{\star}\cdot W_{i-1}  &
=-\sigma(i)\Id_{\mathbb C^i} \text{~for~} 2\leq i\leq n-1,\\
 W_1 \cdot W_1^{\star} & = -\sigma(1);
\end{aligned}
$$
\item there exists a $n\times n$ Hermitian matrix
$H(=W_{n-1}^{\star}\cdot W_{n-1})$ with eigenvalues $$\nu(i)=
-\sum_{j=i}^{n-1}\sigma(j), \forall ~1\leq i \leq n-1 \text{~and~}
\nu(n)=0.$$
\end{enumerate}
\end{lemma}

\begin{proof}
See \cite[Section 3.4]{CB}.
\end{proof}

\begin{proposition}\label{propnou}
Let $\lambda(i)=(\lambda_1(i), \dots, \lambda_n(i)),~1 \leq i \leq
m$ be $m$ weakly decreasing sequences of $n$ reals. Then
$$
\sigma_{\lambda} \in C(Q, \beta)\Longleftrightarrow
(\lambda(1),\dots, \lambda(m) \in \mathcal K(n,m).
$$
\end{proposition}

\begin{proof}
From Proposition \ref{moment}, we know that $\sigma_{\lambda} \in
C(Q, \beta)$ if and only if there exists $W \in \Rep(Q, \beta)$
satisfying the quiver matrix equations $(\dag).$

The matrix equations coming from the first $n-1$ vertices of the
flag $\mathcal F(i)$ are essentially those from Lemma
\ref{linalg}. So, they are equivalent to the existence of
Hermitian matrices $H(i)$ with eigenvalues
$$(\lambda_1(i)-\lambda_n(i),\dots, \lambda_{n-1}(i)-\lambda_n(i),
0).$$ Let $a_1,\dots, a_{m-3}$ denote the main arrows, i.e., those
connecting the central vertices. Taking into account the matrix
equations coming from the main vertices, we see that
$\sigma_{\lambda} \in C(Q, \beta)$ if and only if there exist
Hermitian matrices $H'(i)$ with spectrum $\lambda(i), 1 \leq i
\leq m$ and $n\times n$ complex matrices $W(a_i)$ such that:
$$
\begin{aligned}
&H'(1)+W(a_1)^* \cdot W(a_1)=H'(2),\\
&W(a_1) \cdot W(a_1)^*+W(a_2)\cdot W(a_2)^*=H'(3),\\
&\cdots\\
&H'(m)+W(a_{m-3})^* \cdot W(a_{m-3})=H'(m-1)
\end{aligned}
$$
When writing the last equation of the system above, we assumed
that $m$ is odd. Of course, if $m$ is even, the last equation
looks like $$H'(m)+W(a_{m-3}) \cdot W(a_{m-3})^*=H'(m-1).$$ To
bring the matrix equations above in a for us convenient form, we
can conjugate (if necessary) the equations by unitary matrices.
Also, note that for any $n\times n$ matrix, say $A,$ we have that
$A \cdot A^*$ and $A^* \cdot A$ are both positive semi-definite
and have the same spectrum. Moreover, any positive semi-definite
Hermitian matrix $B$ can be written as $W\cdot W^*$ or $W^* \cdot
W.$

Thus, we obtain that $\sigma_{\lambda} \in C(Q, \beta)$ if and
only if there exist Hermitian matrices $H(i)$ with spectrum
$\lambda(i), 1 \leq i \leq m$ and positive semi-definite $n\times
n$ matrices $B(i)$ such that:
$$
\begin{aligned}
&H(1)+B(1)=H(2),\\
&B(1)+B(2)=H(3),\\
&\cdots\\
&H(m)+B(m-3)=H(m-1).
\end{aligned}
$$
Solving this system of matrix equations for $B(i),$ we have
$$
B(i-1)=\sum_{j=1}^i (-1)^{j+i}H(j), \forall~2\leq i \leq m-2
$$
together with
$$
B(m-3)=H(m-1)-H(m).
$$
Now, the proof follows.
\end{proof}

\begin{proof}[Proof of Theorem \ref{mainthm}]
$(1)\Longleftrightarrow (2)$ This equivalence follows from
Proposition \ref{mainpropnou} and Proposition \ref{propnou}.

$(1)\Longleftrightarrow (3)$ Using Lemma \ref{wtsigmalambda} and
Proposition \ref{propnou} the equivalence follows.

$(3)\Longleftrightarrow (4)$ Note that any long exact sequence
breaks into short exact sequences by taking cokernels. Thus, $(3)$
is equivalent to the existence of short exact sequences
$$0\to M_1\to M_2\to N_1\to 0 ,$$ $$0\to N_1\to M_3\to N_2\to 0,$$
$$\dots$$ $$0\to N_{m-3}\to M_{m-1}\to M_m\to 0,$$ where
$\mu(1),\cdots,\mu(m-3)$ are some partitions of length at most $n$
and $N_1, \dots, N_{m-3}$ are finite abelian $p$-groups of types
$\mu(1), \dots, \mu(m-3).$ This is equivalent to $(4)$ by Klein's
Theorem (see \cite{Kle}).
\end{proof}

\begin{remark}
By definition, we know that $(\lambda(1), \dots, \lambda(m)) \in
\mathcal K(n,m)$ if and only if there exist Hermitian matrices
with prescribed eigenvalues and such that they satisfy a system of
matrix (in)equalities. In principle, one can use the eigenvalue
and the majorization problems (see \cite{CC3} or \cite{F2}) to
find necessary and sufficient Horn type inequalities for each of
the matrix (in)equality defining the cone $\mathcal K(n, m).$ As
we shall see, when we put together these inequalities we obtain a
list of necessary but not sufficient Horn type inequalities. Let
us look at these inequalities when $m=4$ and $n=2.$ In this case,
we want to find inequalities in the parts of $\lambda(1),$
$\lambda(2),$ $\lambda(3),$ $ \lambda(4)$ such that there exist $2
\times 2$ Hermitian matrices $H(1),$ $H(2),$ $H(3),$ $H(4)$ with
eigenvalues $\lambda(1),$ $\lambda(2),$ $\lambda(3),$ $\lambda(4)$
and
$$
H(2)+H(4)=H(1)+H(3)
$$
and
$$
H(1) \leq H(2).
$$
The two conditions above imply the following list of necessary
Horn type inequalities:
$$
|\lambda(2)|+|\lambda(4)|=|\lambda(1)|+|\lambda(3)|,
$$
$$
\lambda_2(2)+\lambda_1(4) \leq \lambda_1(1)+\lambda_1(3),
$$
$$
\lambda_1(2)+\lambda_2(4) \leq \lambda_1(1)+\lambda_1(3),
$$
and
$$
\lambda_2(2)+\lambda_2(4)\leq \lambda_2(1)+\lambda_1(3),
$$
$$
\lambda_2(2)+\lambda_2(4)\leq \lambda_1(1)+\lambda_2(3).
$$
and
$$
\lambda_1(1)\leq \lambda_1(2),~~ \lambda_2(1) \leq \lambda_2(2).
$$

Comparing this list with the one worked out in Example \ref{ex},
we see that the eigenvalue and the majorization problems give
necessary Horn type inequalities which are not sufficient. For
example, take $\lambda(1)=(2,1),$ $\lambda(2)=(3,1),$
$\lambda(3)=(4,1),$ and $\lambda(4)=(2,2).$
\end{remark}

\begin{proof}[Proof of Proposition \ref{genKlycone}]
$(1)$ The chamber inequalities of Lemma \ref{maincoro}{(1)} and
Proposition \ref{propnou} show that
$$
\begin{aligned}
\mathcal K(n,m) & \longrightarrow C(Q, \beta)\times \RR^2 \\
\lambda = (\lambda(1), \dots, \lambda(m)) & \longrightarrow
(\sigma_{\lambda}, \lambda_n(1), \lambda_n(m))
\end{aligned}
$$
is an isomorphism of cones. Since $\beta$ is a Schur root, the
dimension of the cone $C(Q, \beta)$ is the number of the vertices
of the generalized flag quiver minus one and so $(1)$ follows.

$(2)$ This is a consequence of Proposition \ref{mainpropnou}.

\end{proof}

%%%%%%%%%%%%%%%%%%%%%%%%%%%%%%%%%%%%%%%%%%%%%%%%%%%%%%%%%%%%%%%%%%%%%%%%
%%%%%%%%%%%%%%%%%%%%%%%%%%%%%%%%%%%%%%%%%%%%%%%%%%%%%%%%%%%%%%%%%%%%%%%%
\section{Representation theoretic interpretations} \label{reptheoretic}
In this section, we give two representation theoretic
interpretations of the generalized Littlewood-Richardson
coefficients.

\subsection{Parabolic Kazhdan-Lusztig polynomials}
In \cite{LM}, Leclerc and Miyachi obtained some remarkable closed
formulas for certain vectors of the canonical bases of the Fock
space representation of the quantum affine algebra
$U_q(\widehat{{\rm sl}}_n).$ As a direct consequence, they derived
a combinatorial description of certain parabolic affine
Kazhdan-Lusztig polynomials. To state some of their results, we
need to review some definitions from \cite[Section 5]{LM}. Let $v$
be an indeterminate. We denote by $K=\CC(v)$ the field of rational
functions in $v$ and let \textrm{Sym} be the algebra over $K$ of
symmetric functions in a countable set $X$ of variables. Let
$\mathcal P$ be the set of all partitions and $S_{\lambda}$ be the
Schur function labelled by $\lambda \in \mathcal P.$ It is well
known that the functions $S_{\lambda}$ form a linear basis for
\textrm{Sym}. We denote by $\langle \cdot, \cdot \rangle$ the
scalar product for which this basis is orthonormal.

Now, let $N\geq 1$ be an integer and let $A_0,\dots,A_{N-1}$ be
$N$ countable sets of indeterminates. Let
$$
\mathcal S =Sym(A_0, \dots, A_{N-1})
$$
be the algebra over $K$ of functions symmetric in each set
$A_0,\dots, A_{N-1}$ separately. If $\underline{\lambda}=
(\lambda^0,\dots,\lambda^{N-1})\in \mathcal{P}^N,$ consider
$$
S_{\underline{\lambda}}=S_{\lambda^0}(A_0)\cdots
S_{\lambda^{N-1}}(A_{N-1}).
$$

Then $\{ S_{\underline{\lambda}} \mid \underline{\lambda} \in
\mathcal{P}^N \}$ forms a linear basis which is orthonormal with
respect with the induced scalar product. In \cite[Section
5.6]{LM}, the authors introduced a canonical basis $\{
\eta_{\underline{\lambda}}(v) \mid \underline{\lambda} \in
\mathcal{P}^N \}$ and showed that:

\begin{lemma}\cite[Lemma 4]{LM}\label{scalarprod}
For $\underline{\lambda}, \underline{\mu} \in \mathcal{P}^N,$ we
have
$$
\langle S_{\underline{\lambda}}, \eta_{\underline{\mu}}(v) \rangle
= (-v)^{\delta(\underline{\lambda},\underline{\mu})}\sum \prod_{0
\leq j \leq N-1} c_{\alpha^j, \beta^j}^{\mu^j} \cdot c_{\beta^j,
(\alpha^{j+1})'}^{\lambda^j}
$$
where the sum runs through all $\alpha^0,\dots, \alpha^N,$
$\beta^0, \dots, \beta^{N-1}$ in $\mathcal P$ subject to:
$$
|\alpha^i|=\sum_{0 \leq j \leq
i-1}|\lambda^j|-|\mu^j|,\hspace{15pt} |\beta^i|=|\mu^i|+\sum_{0
\leq j \leq i-1}|\mu^j|-|\lambda^j|,
$$
and
$$
\delta(\underline{\lambda},\underline{\mu})=\sum_{0 \leq j \leq
N-2}(N-1-j)(|\lambda^j|-|\mu^j|).
$$
\end{lemma}
Here the convention is that an empty sum is equal to zero. Hence,
$\alpha^0$ is the empty partition, $|\beta^0|=|\mu^0|$ and so $$
c_{\alpha^0, \beta^0}^{\mu^0} \cdot
c_{\beta^0,(\alpha^1)'}^{\lambda^0}=
c_{\mu^0,(\alpha^1)'}^{\lambda^0}.$$ By convention, $\alpha^N$ is
the empty partition and hence $$c_{\alpha^{N-1},
\beta^{N-1}}^{\mu^{N-1}}\cdot
c_{\beta^{N-1},(\alpha^{N})'}^{\lambda^{N-1}}=c_{\alpha^{N-1},\lambda^{N-1}}^{\mu^{N-1}}.$$

Now, let us rewrite the above scalar product using our generalized
Littlewood-Richardson coefficients. It is easy to see that for
$\underline{\lambda}, \underline{\mu} \in \mathcal{P}^N$ we have
$$
\langle S_{\underline{\lambda}}, \eta_{\underline{\mu}}(v) \rangle
= (-v)^{\delta(\underline{\lambda},\underline{\mu})}\cdot
f(\mu^0,\lambda^0,(\mu^1)',(\lambda^1)',\dots, \mu^{N-1},
\lambda^{N-1}).
$$
Note that in the above formula we assumed that $N$ is odd. For $N$
even, just replace $\mu^{N-1}$ and $\lambda^{N-1}$ in $f$ with
$(\mu^{N-1})'$ and $(\lambda^{N-1})'$ respectively.

Next, we explain how these formulas are related to some parabolic
Kazhdan-Lusztig polynomials. Let $w\geq 1$ be an integer and let
$\rho=(\rho_1, \dots, \rho_l)$ be the large $N$-core associated
with $w.$ By $\mathcal P(\rho),$ we denote the set of partitions
with $N$-core $\rho.$ Let $\mathcal{P}(\rho,w)\subseteq \mathcal
P(\rho)$ be the subset of partitions with $N$-weight $\leq w.$ To
each $\lambda \in \mathcal P(\rho),$ one can associate its
$N$-quotient denoted by $\underline{\lambda}=(\lambda^0, \dots,
\lambda^{N-1}).$ For all these definitions, we refer to
\cite[Section 6]{LM}.

\begin{corollary}\cite[Corollary 10]{LM}
Let $\lambda, \mu \in \mathcal{P}(\rho,w).$ Then
$$
d_{\lambda,
\mu}(v)=(-1)^{\delta(\underline{\lambda},\underline{\mu})}\langle
S_{\underline{\lambda}}, \eta_{\underline{\mu}}(v) \rangle \in
\NN[v]
$$
is a parabolic Kazhdan-Lusztig polynomial.
\end{corollary}

Consequently, in this case the coefficient of the Kazhdan-Lusztig
monomial $d_{\lambda, \mu}(v)$ is a generalized
Littlewood-Richardson coefficient. Furthermore, one has that
$d_{\lambda, \mu}(1)$ is a decomposition number of a $q$-Schur
algebra at $q=\sqrt[N]{-1}$ (see also \cite[Theorem 2]{JLM}). Note
that in this case, $d_{\lambda, \mu}(1)$ is a generalized
Littlewood-Richardson coefficient.

In a future paper, we plan to further investigate the connection
between the generalized Littlewood-Richardson coefficients and
decomposition numbers.

\subsection{Multiplicities in representation spaces}
We show that the generalized Littlewood-Richardson coefficients
can be viewed as multiplicities of some irreducible
representations of a product of general linear groups in the
affine coordinate ring of some representation space. For this, let
us consider the alternating type $\mathbb A_m $ quiver with
vertices $1,2,\dots,m$ such that $1$ is a source, $2$ is a sink,
and so on. For example, if $m$ is odd the alternating quiver looks
like:

\begin{diagram}
 1 & \rTo & 2& \cdots & m-1 & \lTo  & m.
\end{diagram}

Now, let $\alpha$ be the dimension vector $\alpha=(n,\dots,n).$
For simplicity, let us write $V(i)=\CC^{n}.$ Without loss of
generality, let us assume that $m$ is odd. Using the
Littlewood-Richardson rule, we can decompose $\CC[\Rep(Q,\alpha)]$
as follows:
$$
\bigoplus
f(\lambda(1),\dots,\lambda(m))\left(S^{\lambda(1)}V(1)\otimes
S^{\lambda(2)}V^{*}(2)\otimes \dots \otimes S^{\lambda(m)}V(m)
\right),
$$
where the sum is taken over all partitions $\lambda(i),~1 \leq i
\leq m$ of length at most $n.$ Thus,
$f(\lambda(1),\dots,\lambda(m))$ is equal to the multiplicity:
$$
\mult_{\GL(\alpha)} \left( S^{\lambda(1)}V(1)\otimes
S^{\lambda(2)}V^{*}(2)\otimes \dots \otimes S^{\lambda(m)}V(m),
\CC[\Rep(Q,\alpha)] \right).
$$

If $m$ is even then $f(\lambda(1),\dots,\lambda(m))$ is equal to
the multiplicity:
$$
\mult_{\GL(\alpha)}\left( S^{\lambda(1)}V(1)\otimes
S^{\lambda(2)}V^{*}(2)\otimes \dots \otimes
S^{\lambda(m)}V^{*}(m), \CC[\Rep(Q,\alpha)] \right).
$$

%We should also point out that for certain partitions the
%non-vanishing of these generalized Littlewood-Richardson numbers
%is equivalent to the non-vanishing of the leading coefficient of
%some parabolic Kazhadan-Lusztig monomials (see Lemma 10.)

\section{Final Remarks}\label{finalrmks}
First, we would like to make some comments regarding the
minimality of our list of Horn type inequalities. When $m=3,$ the
list of necessary and sufficient inequalities from Proposition
\ref{genKlycone}{(2)} is known to be minimal. For $m\geq 4$, this
list of inequalities is not minimal in the sense that it contains
some redundant inequalities. From Remark \ref{facetsrmk} it
follows that the problem concerning the redundancy of our list of
Horn type inequalities comes down to solving the following two
problems.
\newline

\textbf{Problem 1.} Find those $m$-tuples $I = (I_1, \dots, I_m)
\in \mathcal S$ for which the corresponding dimension vectors
$\beta_I$ and $\beta - \beta_I$ are Schur roots.
\newline

If $m=4$, $n=2$, \ \ $I_1 = I_2 = \emptyset,$ ~$I_3 = I_4 =
\{1,2\}$ then  the corresponding dimension vector
$$\beta_I =
\begin{matrix}
 & 0 & 2 \\
0 & 0 & 1 & 1
\end{matrix} $$
is not a Schur root and the redundant inequality is
$|\lambda(4)|\leq |\lambda(3)|$ (see Example \ref{ex}). Now, let
us give examples of tuples $I$ such that both $\beta_I$ and $\beta
- \beta_I$ are Schur roots. Let $I = (I_1, \dots, I_m) \in
\mathcal S(n, m)$ and $|I_i| \geq 1$ for all $1 \leq i \leq m.$ In
case $m
> 3,$ let us assume that
$$
\min \{ |I_3|, |I_{m-2}|, n - |I_3|, n - |I_{m-2}| \} \geq 2,
$$
together with
$$ |I_i| + 1 \leq |I_{i-1}| + |I_{i+1}| \leq n+|I_i|-1,
$$
if $m > 4$ and $3 \leq i \leq m-2.$ Then $\beta_I$ and $\beta -
\beta_I$ are Schur roots. Indeed, we can shrink the generalized
flag quiver so that the restriction of $\beta_I$ to the shrunk
quiver is increasing with jumps all equal to one along the flags.
In this situation, the shrunk dimension vector is indivisible and
lies in the fundamental region. Therefore, it must be a Schur
root. Similarly, one can show that $\beta - \beta_I$ is a Schur
root as well.

The second problem is in fact a particular case of a conjecture in
\cite{DW2}.
\newline

\textbf{Problem 2.} Given $m\geq 3$ partitions $\lambda(1), \dots,
\lambda(m)$, prove that
$$f(\lambda(1), \dots, \lambda(m))=1 \Longleftrightarrow
f(r\lambda(1), \dots, r\lambda(m))=1,$$ for all $r\geq 1.$
\newline

For $m=3,$ this was conjectured by Fulton \cite{F1} and proved by
Knutson, Tao and Woodward \cite{KTW} using puzzles. For arbitrary
$m\geq 3$, the "if" implication is clear. Indeed, we have seen
that $f(r\lambda(1), \dots,
r\lambda(m))=\dim\SI(Q,\beta)_{r\sigma}, \forall ~r\geq 1 $ where
$(Q, \beta)$ is the generalized flag quiver setting. It is easy to
see that $\dim\SI(Q,\beta)_{\sigma}\leq
\dim\SI(Q,\beta)_{r\sigma}, \forall~ r\geq 1$ and so the "if"
implication follows.

%%%%%%%%%%%%%%%%%%%%%%%%%%%%%%%%%%%%%%%%%%%%%%%%%%%%%%%%%%%%%%%%%%%%
%%%%%%%%%%%%%%%%%%%%%%%%%%%%%%%%%%%%%%%%%%%%%%%%%%%%%%%%%%%%%%%%%%%%
\section*{Acknowledgment} I would like to thank William Fulton and
Jerzy Weyman for helpful comments on a preliminary version of this
work. I am especially grateful to my advisor, Harm Derksen, for
drawing my attention to these problems, for his continuing support
and many enlightening discussions on the subject.

%\bibliography{biblio}
\end{document}